\documentclass[aos,rotating,nameyear,dvips]{arximspdf}
\usepackage{dcolumn}

\doi{10.1214/09-AOS781}
\volume{38}
\issue{4}
\pubyear{2010}
\firstpage{2282}
\lastpage{2313}

\makeatletter

\newcolumntype{d}[1]{D{.}{.}{#1}}

\newproclaim{Example}{Example}

\newtheorem{theorem}{Theorem}
\newtheorem{lemma}{Lemma}
\newtheorem{corollary}{Corollary}

\makeatother

\begin{document}
\begin{frontmatter}

\title{Variable selection in nonparametric additive~models}
\runtitle{Nonparametric component selection}

\begin{aug}
\author[A]{\fnms{Jian} \snm{Huang}\thanksref{t1}\ead[label=e1]{jian-huang@uiowa.edu}\corref{}},
\author[B]{\fnms{Joel L.} \snm{Horowitz}\thanksref{t2}\ead[label=e2]{joel-horowitz@northwestern.edu}} and
\author[C]{\fnms{Fengrong} \snm{Wei}\ead[label=e3]{fwei@westga.edu}}
\runauthor{J. Huang, J. L. Horowitz and F. Wei}
\affiliation{University of Iowa, Northwestern University and University
of West Georgia}
\address[A]{J. Huang\\
Department of Statistics\\
\quad and Actuarial Science, 241 SH \\
University of Iowa \\
Iowa City, Iowa 52242\\
USA\\
\printead{e1}} 
\address[B]{J. L. Horowitz\\
Department of Economics\\
Northwestern University\\
2001 Sheridan Road\\
Evanston, Illinois 60208\\
USA\\
\printead{e2}}
\address[C]{F. Wei\\
Department of Mathematics \\
University of West Georgia \\
Carrollton, Georgia 30118 \\
USA\\
\printead{e3}}
\end{aug}

\thankstext{t1}{Supported in part by NIH Grant CA120988 and NSF
Grant DMS-08-05670.}

\thankstext{t2}{Supported in part by NSF Grant SES-0817552.}

\received{\smonth{8} \syear{2009}}
\revised{\smonth{12} \syear{2009}}

%
\begin{abstract}
We consider a nonparametric additive model of a conditional mean
function in which the number of variables and additive components may
be larger than the sample size but the number of nonzero additive
components is ``small'' relative to the sample size. The statistical
problem is to determine which additive components are nonzero. The
additive components are approximated by truncated series expansions
with B-spline bases. With this approximation, the problem of component
selection becomes that of selecting the groups of coefficients in the
expansion. We apply the adaptive group Lasso to select nonzero
components, using the group Lasso to obtain an initial estimator and
reduce the dimension of the problem. We give conditions under which
the group Lasso selects a model whose number of components is
comparable with the underlying model, and the adaptive group Lasso
selects the nonzero components correctly with probability approaching
one as the sample size increases and achieves the optimal rate of
convergence. The results of Monte Carlo experiments show that the
adaptive group Lasso procedure works well with samples of moderate
size. A data example is used to illustrate the application of the
proposed method.
\end{abstract}

%
\begin{keyword}[class=AMS]
\kwd[Primary ]{62G08}
\kwd{62G20}
\kwd[; secondary ]{62G99}.
\end{keyword}
\begin{keyword}
\kwd{Adaptive group Lasso}
\kwd{component selection}
\kwd{high-dimensional data}
\kwd{nonparametric regression}
\kwd{selection consistency}.
\end{keyword}

\end{frontmatter}

\section{Introduction}

Let $(Y_i, {\mathbf X}_i), i=1, \ldots, n$, be random vectors that are
independently
and identically distributed as $(Y, {\mathbf X})$, where $Y$ is a
response variable and ${\mathbf X}=(X_{1}, \ldots, X_{p})'$ is a
$p$-dimensional covariate vector. Consider the nonparametric
additive model
%
%
\begin{equation}
\label{GamA} Y_i= \mu+ \sum_{j=1}^{p} f_j(X_{ij}) + \varepsilon_i,
\end{equation}
where $\mu$ is an intercept term, $X_{ij}$ is the $j$th component of
$X_{i}$, the $f_j$'s are unknown
functions, and $\varepsilon_i$ is an unobserved random variable
with mean zero and finite variance $\sigma^2$. Suppose that some of the
additive components $f_j$ are zero. The problem addressed in this paper
is to distinguish the nonzero components from the zero components and
estimate the nonzero components. We allow the possibility that $p$ is
larger than
the sample size $n$, which we represent by letting $p$ increase as $n$
increases. We propose a penalized method for variable selection in
(\ref
{GamA}) and show that the proposed method can correctly select
the nonzero components with high probability.

There has been much work on penalized methods for variable
selection and estimation with high-dimensional data. Methods that have
been proposed include the bridge estimator [Frank and Friedman
(\citeyear{FF93}), Huang, Horowitz and Ma (\citeyear{HHM08})]; least absolute
shrinkage and selection operator or Lasso [Tibshirani (\citeyear{Tibshirani96})], the
smoothly clipped absolute deviation (SCAD) penalty [Fan and Li
(\citeyear{FL01}), Fan and Peng (\citeyear{FP04})],
and the minimum concave penalty [Zhang (\citeyear{Zhang07})].
Much progress has been made in understanding the statistical properties
of these methods. In particular, many authors have studied the variable
selection, estimation and prediction properties of the Lasso in
high-dimensional settings. See, for example, Meinshausen and
B\"{u}hlmann (\citeyear{MB06}), Zhao and Yu (\citeyear{ZhaoYu06}), Zou (\citeyear{Zou06}),
Bunea, Tsybakov and Wegkamp (\citeyear{BTW07}),
Meinshausen and Yu (\citeyear{MY09}), Huang, Ma and Zhang
(\citeyear{HMZ08}), van de Geer (\citeyear{vandeGeer08})
and Zhang and Huang (\citeyear{ZH08}), among others. All these
authors assume a linear or other parametric model. In many
applications, however, there is little a priori justification for
assuming that the effects of covariates take a linear form or belong to
any other known, finite-dimensional parametric family. For example, in
studies of economic development, the effects of covariates on the
growth of gross domestic product can be nonlinear. Similarly, there is
evidence of nonlinearity in the gene expression data used in the
empirical example in Section~\ref{sec5}.\looseness=1

There is a large body of literature on estimation in nonparametric
additive models. For example,
Stone (\citeyear{Stone85}, \citeyear{Stone86}) showed that additive
spline estimators achieve the same optimal rate of convergence for a
general fixed $p$ as for $p=1$. Horowitz and Mammen (\citeyear{HM04})
and Horowitz, Klemel\"{a} and Mammen (\citeyear{HKM06}) showed that if $p$ is fixed
and mild regularity conditions hold, then oracle-efficient estimates of
the $f_j$'s can be obtained by a two-step procedure. Here, oracle
efficiency means that the estimator of each $f_j$ has the same
asymptotic distribution that it would have if all the other $f_j$'s
were known. However, these papers do not discuss variable selection in
nonparametric additive models.

Antoniadis and Fan (\citeyear{AnFan01}) proposed a group SCAD approach for
regularization in wavelets approximation.
Zhang et al. (\citeyear{Zhangetal04}) and Lin and Zhang (\citeyear{LinZhang06})
have investigated the use
of penalization methods in smoothing spline ANOVA with a fixed number
of covariates. Zhang et al. (\citeyear{Zhangetal04}) used a Lasso-type penalty but did
not investigate model-selection consistency. Lin and Zhang (\citeyear{LinZhang06})
proposed the component selection and smoothing operator (COSSO) method
for model selection and estimation in multivariate nonparametric
regression models. For fixed $p$, they showed that the COSSO estimator
in the additive model converges at the rate $n^{-d/(2d+1)}$, where $d$
is the order of smoothness of the components. They also showed that, in
the special case of a tensor product design, the COSSO correctly
selects the nonzero additive components with high probability. Zhang
and Lin (\citeyear{ZhangLin06}) considered the COSSO for nonparametric regression in
exponential families.

Meier, van de Geer and B\"{u}hlmann (\citeyear{MGB09}) treat variable selection in
a nonparametric additive model in which the numbers of zero and nonzero
$f_j$'s may both be larger than $n$. They propose a penalized
least-squares estimator for variable selection and estimation. They
give conditions under which, with probability approaching 1, their
procedure selects a set of $f_j$'s containing all the additive
components whose distance from zero in a certain metric exceeds a
specified threshold. However, they do not establish model-selection
consistency of their procedure. Even asymptotically, the selected set
may be larger than the set of nonzero $f_j$'s. Moreover, they impose a
compatibility condition that relates the levels and smoothness of the
$f_j$'s. The compatibility condition does not have a straightforward,
intuitive interpretation and, as they point out, cannot be checked
empirically. Ravikumar et al. (\citeyear{Ravikumaretal09}) proposed a
penalized approach for variable selection in nonparametric additive
models. In their approach, the penalty is imposed on the $\ell_2$ norm
of the nonparametric components, as well as the mean value of the
components to ensure identifiability.
In their theoretical results, they require that the eigenvalues of a
``design matrix'' be bounded away from zero and infinity, where the
``design matrix'' is formed from the basis functions for the nonzero
components. It is not clear whether this condition holds in general,
especially when the number of nonzero components diverges with $n$.
Another critical condition required in the results of Ravikumar et al.
(\citeyear{Ravikumaretal09}) is similar to the irrepresentable condition of Zhao and Yu
(\citeyear{ZhaoYu06}). It is not clear for what type of basis functions this condition
is satisfied. We do not require such a condition in our results on
selection consistency of the adaptive group Lasso.

Several other recent papers have also considered variable selection in
nonparametric models. For example, Wang, Chen and Li (\citeyear{WCL07}) and
Wang and Xia (\citeyear{WangXia08}) considered the use of group Lasso and SCAD methods
for model selection and estimation in varying coefficient models with a
fixed number of coefficients and covariates. Bach (\citeyear{Bach07}) applies what
amounts to the group Lasso to a nonparametric
additive model with a fixed number of covariates. He established model
selection consistency under conditions that are considerably more
complicated than the ones we require for a possibly diverging number of
covariates.

In this paper, we propose to use the adaptive group Lasso for variable
selection in (\ref{GamA}) based on a spline approximation to the
nonparametric components. With this approximation, each nonparametric
component is represented by a linear combination of spline basis
functions. Consequently, the problem of component selection becomes
that of selecting the groups of coefficients in the linear
combinations. It is natural to apply the group Lasso method, since it
is desirable to take into the grouping structure in the approximating
model. To achieve model selection consistency, we apply the group Lasso
iteratively as follows.
First, we use the group Lasso to obtain an initial estimator and reduce
the dimension of the problem. Then we use the adaptive group Lasso to select
the final set of nonparametric components.
The adaptive group Lasso is a simple generalization of the adaptive Lasso
[Zou (\citeyear{Zou06})] to the method of the group Lasso [Yuan and Lin (\citeyear{YL06})].
However, here we apply this approach to nonparametric additive modeling.

We assume that the number of nonzero $f_j$'s is fixed. This enables us
to achieve model selection consistency under simple assumptions that
are easy to interpret. We do not have to impose compatibility or
irrepresentable conditions, nor do we need to assume conditions on the
eigenvalues of certain matrices formed from the spline basis functions.
We show that the group Lasso selects a model whose number of components
is bounded with probability approaching one by a constant that is
independent of the sample size. Then using the group Lasso result as
the initial estimator, the adaptive group Lasso selects the correct
model with probability approaching 1 and achieves the optimal rate of
convergence for nonparametric estimation of an additive model.

The remainder of the paper is organized as follows. Section \ref{sec2}
describes
the group Lasso and the adaptive group Lasso for variable selection in
nonparametric additive models. Section \ref{sec3} presents the asymptotic
properties of these methods in ``large $p$, small $n$'' settings.
Section \ref{sec4} presents the results of simulation studies to
evaluate the
finite-sample performance of these methods. Section \ref{sec5}
provides an
illustrative application, and Section \ref{sec6} includes concluding remarks.
Proofs of the results stated in Section \ref{sec3} are given in
the \hyperref[app]{Appendix}.

\section{Adaptive group Lasso in nonparametric additive models}
\label{sec2}

We describe a two-step approach that uses the group Lasso for variable
selection based on a spline representation of each component in
additive models. In the first step, we use the standard group Lasso to
achieve an initial reduction of the dimension in the model and obtain
an initial estimator of the nonparametric components. In the second
step, we use the adaptive group Lasso to achieve consistent selection.

Suppose that each $X_j$ takes values in $[a, b]$ where $a < b$ are
finite numbers. To ensure unique identification of the $f_j$'s, we
assume that $\mathrm{E}f_j(X_j) =0, 1 \le j \le p$.
Let $a=\xi_0 < \xi_1 < \cdots< \xi_K < \xi_{K+1}=b$
be a partition of $[a, b]$ into $K$ subintervals
$I_{Kt} = [\xi_t, \xi_{t+1}), t=0, \ldots, K-1$, and
$I_{KK}=[\xi_K, \xi_{K+1}]$, where $K \equiv K_n = n^v$ with $0 < v < 0.5$
is a positive integer such that
$\max_{1 \le k \le K+1}|\xi_k-\xi_{k-1}| = O(n^{-v})$.
Let $\mathcal{S}_n$ be the space of polynomial splines of degree $l
\ge1$
consisting of functions $s$ satisfying:
(i) the restriction of $s$ to $I_{Kt}$ is a polynomial of degree $l$
for $1 \le t \le K$; (ii) for $ l \ge2$ and $ 0 \le l' \le l-2$,
$s$ is $l'$ times continuously differentiable on $[a, b]$.
This definition is phrased after Stone (\citeyear{Stone85}), which is
a descriptive version of Schumaker (\citeyear{Schumaker81}), page 108, Definition 4.1.

There exists a normalized B-spline basis $\{\phi_{k}, 1 \le k \le m_n
\}$ for $\mathcal{S}_n$, where $m_n \equiv K_n + l$ [Schumaker (\citeyear{Schumaker81})].
Thus, for
any $f_{nj} \in\mathcal{S}_n$, we can write
%
%
\begin{equation}
\label{RepA}
f_{nj}(x) = \sum_{k=1}^{m_n} \beta_{jk}\phi_{k}(x),\qquad
1 \le j \le p.
\end{equation}
Under suitable smoothness assumptions, the $f_{j}$'s can be well
approximated by functions in $\mathcal{S}_n$. Accordingly, the variable
selection method described in this paper is based on the representation
(\ref{RepA}).

Let $\|{\mathbf a}\|_2 \equiv(\sum_{j=1}^{m}|a_j|^{2} )^{1/2}$ denote
the $\ell_2$ norm of any vector ${\mathbf a}\in\mathbb{R}^m$.
Let $\bolds\beta_{nj}=(\beta_{j1}, \ldots, \beta_{jm_n})'$ and
$\bolds\beta_n=(\bolds\beta_{n1}', \ldots, \bolds\beta
_{np}')'$. Let $w_n= (w_{n1},
\ldots, w_{np})'$ be a given vector of
weights, where $0 \le w_{nj}\le\infty, 1\le j \le p$. Consider the
penalized least squares criterion
%
%
\begin{equation}\quad
\label{CA} L_n(\mu, \bolds\beta_n) =
\sum_{i=1}^n \Biggl[Y_i -\mu-\sum_{j=1}^{p} \sum_{k=1}^{m_n} \beta
_{jk}\phi_{k}(X_{ij}) \Biggr]^2 + \lambda_n \sum_{j=1}^{p} w_{nj}
\|\bolds\beta_{nj}\|_2,
\end{equation}
where $\lambda_n$ is a penalty parameter. We study the estimators that
minimize $L_n(\mu, \bolds\beta_n)$ subject to the constraints
%
%
\begin{equation}
\label{RA} 
\sum_{i=1}^n \sum_{k=1}^{m_n}\beta_{jk}\phi_{k}(X_{ij})=0,\qquad 1 \le j
\le p.
\end{equation}
These centering constraints are sample analogs of the identifying
restriction $\mathrm{E}f_j(X_j) = 0, 1 \le j \le p$. We can convert (\ref
{CA}) and (\ref{RA}) to an unconstrained optimization problem by centering
the response and the basis functions. Let
%
%
\begin{equation}
\label{CenterA}
\bar{\phi}_{jk} = \frac{1}{n}\sum_{i=1}^n \phi_k(X_{ij}),\qquad
\psi_{jk}(x) =\phi_k(x)-\bar{\phi}_{jk}.
\end{equation}
For simplicity and without causing confusion, we simply write
$\psi_{k}(x)=\psi_{jk}(x)$.
Define
\[
Z_{ij}= (\psi_{1}(X_{ij}), \ldots,
\psi_{m_n}(X_{ij}) )'.
\]
So, $Z_{ij}$ consists of values of the (centered) basis functions at
the $i$th observation of the $j$th covariate. Let $ {\mathbf Z}_j = (Z_{1j},
\ldots, Z_{nj})'$ be the $n \times m_n$ ``design'' matrix corresponding
to the $j$th covariate. The total ``design'' matrix is ${\mathbf
Z}=({\mathbf Z}_1,
\ldots
, {\mathbf Z}_p)$. Let ${\mathbf Y}=(Y_1-\overline{Y}, \ldots,
Y_n-\overline{Y})'$. With this
notation, we
can write
%
%
\begin{equation}
\label{CB} L_n(\bolds\beta_n;\lambda) = \|{\mathbf Y}-{\mathbf
Z}\bolds\beta_n\|_2^2 +
\lambda_n
\sum_{j=1}^{p} w_{nj} \|\bolds\beta_{nj}\|_2.
\end{equation}
Here, we have dropped $\mu$ in the argument of $L_n$. With the centering,
$\widehat{\mu}= \overline{Y}$. Then minimizing (\ref{CA}) subject
to (\ref{RA}) is
equivalent to minimizing (\ref{CB}) with respect to $\bolds\beta
_n$, but the
centering constraints are not needed for (\ref{CB}).

We now describe the two-step approach to component selection in
the nonparametric additive model (\ref{GamA}).

\textit{Step} 1. Compute the group Lasso estimator. Let
\[
L_{n1}(\bolds\beta_n,\lambda_{n1}) = \|{\mathbf Y}-{\mathbf
Z}\bolds\beta_n\|_2^2 + \lambda_{n1}
\sum_{j=1}^{p} \|\bolds\beta_{nj}\|_2.
\]
This objective function is the special case of (\ref{CB}) that is
obtained by setting $w_{nj}=1$, $1\le j\le p$. The group Lasso estimator is
\(
\widetilde{\bolds\beta}_n \equiv\widetilde{\bolds\beta
}_n(\lambda_{n1}) = \arg\min_{\bolds\beta_n}
L_{n1}(\bolds\beta_n;\lambda_{n1}).
\)

\textit{Step} 2. Use the group Lasso estimator $\widetilde{\bolds
\beta}_n$ to obtain the
weights by setting
\[
w_{nj}= \cases{
\|\widetilde{\bolds\beta}_{nj}\|_2^{-1}, &\quad if $\|\widetilde{\bolds
\beta}_{nj}\|_2 > 0$, \cr
\infty, &\quad if $\|\widetilde{\bolds\beta}_{nj}\|_2 =0$.}
\]
The adaptive group Lasso objective function is
\[
L_{n2}(\bolds\beta_n;\lambda_{n2}) = \|{\mathbf Y}-{\mathbf
Z}\bolds\beta_n\|_2^2 + \lambda_{n2}
\sum_{j=1}^{p} w_{nj} \|\bolds\beta_{nj}\|_2.
\]
Here, we define $0 \cdot\infty=0$. Thus, the components not selected
by the group Lasso are not included in Step 2.
The adaptive group Lasso estimator is $\widehat{\bolds\beta}_n
\equiv\widehat{\bolds\beta}
_n(\lambda_{n2})= \arg\min_{\bolds\beta_n} L_{n2}(\bolds\beta
_n;\lambda_{n2})$.
Finally, the adaptive group Lasso estimators of $\mu$ and $f_j$ are
\[
\widehat{\mu}_n = \overline{Y}\equiv n^{-1}\sum_{i=1}^n Y_i,
\widehat{f}_{nj}(x) = \sum_{k=1}^{m_n} \widehat{\beta}_{jk} \psi
_{k}(x),\qquad 1 \le
j \le p.
\]

\section{Main results}\label{sec3}

This section presents our results on the asymptotic properties of the
estimators defined in Steps 1 and 2 of Section \ref{sec2}.

Let $k$ be a nonnegative integer, and let $\alpha\in(0, 1]$ be such
that $d = k+\alpha> 0.5$. Let $\mathcal{F}$ be the class of functions
$f$ on $[0,1]$
whose $k$th derivative $f^{(k)}$ exists and satisfies a
Lipschitz condition of order $\alpha$:
\[
\bigl|f^{(k)}(s)- f^{(k)}(t)\bigr| \le C |s-t|^{\alpha}\qquad
\mbox{for } s, t \in[a, b].
\]

In (\ref{GamA}), without loss of generality, suppose that the
first $q$ components are nonzero, that is, $f_j(x) \neq0, 1 \le
j \le q$, but $f_j(x)\equiv0, q+1 \le j \le p$.
Let $A_1=\{1, \ldots, q\}$ and $A_0=\{q+1, \ldots, p\}$. Define
$\|f\|_2 = [\int_a^bf^2(x)\,dx]^{1/2}$ for any function $f$, whenever the
integral exists.

We make the following assumptions.

(A1) The number of nonzero components $q$ is fixed and
there is a constant $c_f > 0$ such that
${\min_{1\le j\le q}} \|f_j\|_2 \ge c_f$.

(A2) The random variables $\varepsilon_1,\ldots, \varepsilon_n$
are independent and
identically distributed with $\mathrm{E}
\varepsilon_i=0$ and $\operatorname{Var}(\varepsilon_i)=\sigma
^2$. Furthermore, their tail
probabilities satisfy $P(|\varepsilon_i| > x) \le K \exp(-C x^{2}),
i=1,\ldots, n$, for all $x \ge0$ and for constants $C$ and~$K$.

(A3) $\mathrm{E}f_j(X_j)=0$ and $f_j \in\mathcal{F}, j=1, \ldots, q$.

(A4) The covariate vector $X$ has a continuous density
and there exist constants $C_1$ and $C_2$ such that the density
function $g_j$ of $X_j$ satisfies $0 < C_1 \le g_j(x) \le C_2 < \infty$
on $[a, b]$ for every $1 \le j \le p$.

We note that (A1), (A3) and (A4) are standard conditions for
nonparametric additive models. They would be needed to estimate the
nonzero additive components at the optimal $\ell_2$ rate of convergence
on $[a, b]$, even if $q$ were fixed and known. Only (A2) strengthens
the assumptions needed for nonparametric estimation of a nonparametric
additive model.
While condition (A1) is reasonable in most applications, it would be
interesting to relax this condition and investigate the case when the
number of nonzero components can also increase with the sample size.
The only technical reason that we assume this condition is related to
Lemma \ref{LemPB} given in the \hyperref[app]{Appendix}, which is concerned with the properties
of the smallest and largest eigenvalues of the ``design matrix'' formed
from the spline basis functions. If this lemma can be extended to the
case of a divergent number of components, then (A1) can be relaxed.
However, it is clear that there needs to be restriction on
the number of nonzero components to ensure model identification.

\subsection{Estimation consistency of the group Lasso}\label{sec31}

In this section, we consider the selection and estimation properties of
the group Lasso estimator. Define ${\widetilde A}_1=\{j\dvtx\|
\widetilde{\bolds\beta}_{nj}\|_2
\neq
0, 1\le j \le p\}$. Let $|A|$ denote the cardinality of any set $A
\subseteq\{1, \ldots, p\}$.
\begin{theorem}
\label{ThmA} Suppose that \textup{(A1)} to \textup{(A4)} hold and
$\lambda_{n1} \ge C \sqrt{n \log(p m_n)} $ for a sufficiently large
constant $C$.

\begin{longlist}
\item With probability converging to 1, $|{\widetilde A}_1| \le M_1
|A_1|=M_1q$ for
a finite constant $M_1 > 1$.

\item If ${m_n^2 \log(p m_n)}/{n}\rightarrow0$ and $(\lambda_{n1}^2
m_n)/n^2 \rightarrow0$ as $n \rightarrow\infty$, then all the nonzero
$\bolds\beta_{nj}, 1\le j \le q$, are selected with probability converging
to one.\vspace*{15pt}

\item\mbox{}\vspace*{-30pt}
\begin{eqnarray*}
\sum_{j=1}^{p} \|\widetilde{\bolds\beta}_{nj}-\bolds\beta_{nj}\|
_2^2 &=&
O_p \biggl(\frac{m_n^2 \log(p m_n)}{n} \biggr) +
O_p \biggl(\frac{m_n}{n} \biggr)\\
&&{} + O \biggl(\frac{1}{
m_n^{2d-1}} \biggr) + O \biggl(\frac{4 m_n^2\lambda_{n1}^2
}{n^2} \biggr).
\end{eqnarray*}
\end{longlist}
\end{theorem}

Part (i) of Theorem \ref{ThmA} says that, with probability approaching
1, the group Lasso selects a model whose dimension is a constant
multiple of the number of nonzero additive components $f_j$, regardless
of the number of additive components that are zero. Part (ii) implies
that every nonzero coefficient will be selected with high probability.
Part (iii) shows that the difference between the coefficients in the
spline representation
of the nonparametric functions in (\ref{GamA}) and their estimators
converges to zero in probability. The rate of convergence is determined by
four terms: the stochastic error in estimating the nonparametric
components (the first term) and the intercept $\mu$ (the second term),
the spline approximation error (the third term) and the bias due to penalization
(the fourth term).

Let $\widetilde{f}_{nj}(x)=\sum_{j=1}^{m_n}\widetilde{\beta
}_{jk}\psi(x), 1 \le j \le p$.
The following theorem is a consequence of Theorem \ref{ThmA}.
\begin{theorem}
\label{ThmB}
Suppose that \textup{(A1)} to \textup{(A4)} hold and that $\lambda_{n1} \ge\break C\sqrt
{n\log
(pm_n)}$ for a sufficiently large constant $C$.
Then:

\begin{longlist}
\item Let ${\widetilde A}_f=\{j\dvtx\|\widetilde{f}_{nj}\|_2 > 0, 1\le
j\le p\}$.
There is a constant $M_1 > 1$ such that, with probability converging to
1, $|{\widetilde A}_f| \le M_1 q$.

\item If $(m_n\log(pm_n))/n \rightarrow0$ and $(\lambda_{n1}^2 m_n)/n^2
\rightarrow0$ as $n \rightarrow\infty$, then all the nonzero additive
components $f_j, 1\le j \le q$, are selected with probability
converging to one.\vspace*{10pt}

\item\mbox{}\vspace*{-27pt}
\begin{eqnarray*}
\|\widetilde{f}_{nj}-f_j\|_2^2 &=&
O_p \biggl(\frac{ m_n \log(p m_n)}{n} \biggr)+ O_p \biggl(\frac
{1}{n} \biggr) \\
&&{} + O \biggl(\frac{1}{
m_n^{2d}} \biggr) + O \biggl(\frac{4 m_n\lambda_{n1}^2}{n^2} \biggr),
\qquad
j \in{\widetilde A}_2,
\end{eqnarray*}
where ${\widetilde A}_2 = A_1 \cup{\widetilde A}_1$.
\end{longlist}
\end{theorem}

Thus, under the conditions of Theorem \ref{ThmB}, the group Lasso
selects all the nonzero additive components with high probability. Part
(iii) of the theorem gives the rate of convergence of the group Lasso
estimator of the nonparametric components.

For any two sequences $\{a_n, b_n, n=1, 2, \ldots\}$, we write $a_n
\asymp b_n$ if there
are constants $0 < c_1 < c_2 < \infty$ such that $c_1 \le a_n/b_n \le c_2$
for all $n$ sufficiently large.


We now state a useful corollary of Theorem \ref{ThmB}.
\begin{corollary}
\label{ColA}
Suppose that \textup{(A1)} to \textup{(A4)} hold. If $\lambda_{n1} \asymp\sqrt{n \log(pm_n)}$
and $m_n \asymp n^{1/(2d+1)}$, then:

\begin{longlist}
\item If $n^{-2d/(2d+1)}\log(p) \rightarrow0$ as $n \rightarrow\infty$,
then with probability converging to one, all the nonzero components
$f_j, 1\le j \le q$, are selected and the number of selected
components is no more than $M_1 q$.\vspace*{10pt}

\item\mbox{}\vspace*{-22pt}
\[
\|\widetilde{f}_{nj}-f_j\|_2^2 =O_p\bigl(n^{-2d/(2d+1)}\log(p m_n)\bigr),\qquad j \in
{\widetilde A}_2.
\]
\end{longlist}
\end{corollary}

For the $\lambda_{n1}$ and $m_n$ given in Corollary \ref{ColA}, the
number of zero components can be as large as
$\exp(o(n^{2d/(2d+1)}))$. For example, if each $f_j$ has continuous
second derivative ($d=2$), then it is $\exp(o( n^{4/5}))$, which can be
much larger than $n$.

\subsection{Selection consistency of the adaptive group Lasso}
\label{sec32}

We now consider the properties of the adaptive group Lasso.
We first state a general result concerning the selection consistency
of the adaptive group Lasso, assuming an initial consistent estimator
is available. We then apply to the case when the group Lasso is used as the
initial estimator. We make the following assumptions.

(B1) The initial estimators $\widetilde{\bolds\beta}_{nj}$ are
$r_n$-consistent at zero:
\[
r_n \max_{j \in A_0}\|\widetilde{\bolds\beta}_{nj}\|_2 = O_P(1),\qquad
r_n\to\infty,
\]
and there exists a constant $c_b > 0$ such that
\[
\mathrm{P}\Bigl({\min_{j \in A_1}} \|\widetilde{\bolds\beta}_{nj}\|_2 \ge
c_b b_{n1}\Bigr)
\rightarrow1,
\]
where $b_{n1}={\min_{j \in A_1}} \|\bolds\beta_{nj}\|_2$.

(B2) Let $q$ be the number of nonzero components and $s_n=p-q$
be the number of zero components. Suppose that:
\begin{eqnarray*}
\mbox{(a)\hspace*{151pt}} \frac{m_n}{n^{1/2}} +
\frac{\lambda_{n2} m_n^{1/4}}{n } &=& o(1),\hspace*{150pt} \\
\mbox{(b)\hspace*{60pt}}\frac{n^{1/2}\log^{1/2}(s_nm_n)}{\lambda_{n2} r_n}+
\frac{n}{\lambda_{n2} r_nm_n^{(2d+1)/2}} 
&=& o(1).\hspace*{60pt}
\end{eqnarray*}

We state condition (B1) for a general initial estimator, to highlight the
point that the availability of an $r_n$-consistent estimator at zero is
crucial for the adaptive group Lasso to be selection consistent. In
other words, any initial estimator satisfying (B1) will ensure that the
adaptive group Lasso (based on this initial estimator) is selection
consistent, provided that certain regularity conditions are satisfied.
We note that it follows immediately from Theorem \ref{ThmA} that the
group Lasso estimator satisfies (B1). We will come back to this point below.

For $\widehat{\bolds\beta}_n\equiv(\widehat{\bolds\beta}_{n1}',
\ldots, \widehat{\bolds\beta}_{np}')'$
and $\bolds\beta_n\equiv(\bolds\beta_{n1}', \ldots, \bolds
\beta_{np}')'$, we
say $\widehat{\bolds\beta}_n=_0 \bolds\beta_n$ if $\operatorname
{sgn}_0(\|\widehat{\bolds\beta}_{nj}\|) =
\operatorname{sgn}_0(\|\bolds\beta_{nj}\|), 1 \le j \le p$, where
$\operatorname{sgn}_0(|x|)=1$
if $|x| > 0$ and $=0$ if $|x|=0$.
\begin{theorem}
\label{ThmC} Suppose that conditions \textup{(B1)}, \textup{(B2)} and \textup{(A1)--(A4)}
hold. Then:
\begin{eqnarray*}
\mbox{\textup{(i)}\hspace*{101pt}}\mathrm{P}(\widehat{\bolds\beta}_n =_0 \bolds\beta_n)
&\rightarrow&1.\hspace*{100pt}
\\
\mbox{\textup{(ii)}\hspace*{79pt}}\sum_{j=1}^{q} \|\widehat{\bolds\beta}_{nj}-\bolds\beta_{nj}\|
_2^2 &=&
O_p \biggl(\frac{m_n^2 }{n} \biggr) +
O_p \biggl(\frac{m_n}{n} \biggr)\\
&&{}+O \biggl(\frac{1}{
m_n^{2d-1}} \biggr) + O \biggl(\frac{4 m_n^2\lambda_{n2}^2
}{n^2} \biggr).\hspace*{79pt}
\end{eqnarray*}
\end{theorem}

This theorem is concerned with the selection and estimation
properties of the adaptive group Lasso in terms of $\widehat{\bolds
\beta}_n$.
The following theorem states the results in terms of the estimators
of the nonparametric components.
\begin{theorem}
\label{ThmD} Suppose that conditions \textup{(B1)}, \textup{(B2)} and \textup{(A1)--(A4)}
hold. Then:
\[
\mbox{\textup{(i)}\hspace*{63pt}}\mathrm{P}(\|\widehat{f}_{nj}\|_2 > 0, j \in A_1 \mbox{ and }
\|\widehat{f}_{nj}\|_2=0, j \in A_0 ) \rightarrow 1.\hspace*{50pt}
\]
\begin{eqnarray*}
\\[-20pt]
\mbox{\textup{(ii)}\hspace*{75pt}}\sum_{j=1}^{q} \|\widehat{f}_{nj}- f_{j}\|_2^2
&=&
O_p \biggl(\frac{m_n }{n} \biggr) +
O_p \biggl(\frac{1}{n} \biggr)\\[-5pt]
&&{}+O \biggl(\frac{1}{
m_n^{2d}} \biggr) + O \biggl(\frac{4 m_n\lambda_{n2}^2
}{n^2} \biggr).\hspace*{75pt}
\end{eqnarray*}
\end{theorem}

Part (i) of this theorem states that the
adaptive group Lasso can consistently distinguish nonzero
components from zero components. Part (ii) gives an upper
bound on the rate of convergence of the estimator.

We now apply the above results to our proposed procedure described in
Section~\ref{sec2}, in which we first obtain the the group Lasso
estimator and then use it as the initial estimator in the adaptive
group Lasso.

By Theorem \ref{ThmA}, if $\lambda_{n1}\asymp\sqrt{n \log(pm_n)}$ and
$m_n \asymp n^{1/(2d+1)}$ for $ d \ge1$, then the group Lasso estimator
satisfies (B1) with
$r_n \asymp n^{d/(2d+1)}/\sqrt{ \log(pm_n)}$. In this case, (B2)
simplifies to
%
%
\begin{equation}
\label{AdaA}
\frac{\lambda_{n2} }{n^{(8d+3)/(8d+4)}} = o(1) \quad\mbox{and}\quad
\frac{n^{1/(4d+2)}\log^{1/2}(pm_n)}{\lambda_{n2}}=o(1).
\end{equation}

We summarize the above discussion in the following corollary.
\begin{corollary}
\label{ColB}
Let the group Lasso estimator $\widetilde{\bolds\beta}_n\equiv
\widetilde{\bolds\beta}_n(\lambda_{n1})$
with $\lambda_{n1}\asymp\sqrt{n\log(pm_n)}$ and $m_n \asymp
n^{1/(2d+1)}$ be the initial estimator in the adaptive group Lasso.
Suppose that the conditions of Theorem \ref{ThmA} hold.
If $\lambda_{n2}\le O(n^{1/2})$ and satisfies (\ref{AdaA}),
then the adaptive group Lasso consistently selects the nonzero
components in (\ref{GamA}), that is, part \textup{(i)} of Theorem \ref{ThmD}
holds. In addition,
\[
\sum_{j=1}^{q} \|\widehat{f}_{nj}- f_{j}\|_2^2 =
O_p \bigl(n^{-2d/(2d+1)} \bigr).
\]
\end{corollary}

This corollary\vspace*{1pt} follows directly from Theorems \ref{ThmA} and \ref{ThmD}.
The largest $\lambda_{n2}$ allowed is $\lambda_{n2}=O(n^{1/2})$.
With this $\lambda_{n2}$, the first equation in (\ref{CB}) is satisfied.
Substitute it into the second equation in (\ref{CB}), we obtain $p = \exp
(o(n^{2d/(2d+1)}))$, which is the largest $p$ permitted and can be
larger than $n$. Thus, under the conditions of this corollary, our
proposed adaptive group Lasso estimator using the group Lasso as the
initial estimator is selection consistent and achieves optimal rate
of convergence even when $p$ is larger than $n$. Following model
selection, oracle-efficient, asymptotically normal estimators of the
nonzero components can be obtained by using existing methods.

\section{Simulation studies}\label{sec4}

We use simulation to evaluate the performance of the adaptive group
Lasso with regard to variable selection. The generating model is
%
%
\begin{equation}
\label{GModela}
y_{i}=f(x_i)+\varepsilon_i \equiv\sum
_{j=1}^{p}f_{j}(x_{ij})+\varepsilon
_{i},\qquad i=1,\ldots, n.
\end{equation}

Since $p$ can be larger than $n$, we consider two ways to select the
penalty parameter, the BIC
[Schwarz (\citeyear{Schwarz78})] and the EBIC [Chen and Chen (\citeyear{ChCh08}, \citeyear{ChCh09})]. The BIC is
defined as
\[
\mathit{BIC}(\lambda)=\log(\mathrm{RSS}_{\lambda})+ df_{\lambda}\cdot\frac
{\log n}{n}.
\]
Here, $\mathrm{RSS}_{\lambda}$ is the residual sum of squares for a given
$\lambda$, and the degrees of freedom $df_{\lambda} = \hat
{q}_{\lambda}
m_n$, where $\hat{q}_{\lambda}$ is the number of nonzero estimated
components for the given $\lambda$.
The EBIC is defined as
\[
\mathit{EBIC}(\lambda)=\log(\mathrm{RSS}_{\lambda})+df_{\lambda} \cdot\frac{\log
n}{n}+\nu\cdot df_{\lambda} \cdot\frac{\log p}{n},
\]
where $0 \le\nu\le1$ is a constant. We use $\nu=0.5$.

We have also considered two other possible ways of defining df: (a)
using the trace of a linear smoother based on a quadratic
approximation; (b) using the number of estimated nonzero components. We
have decided to use the definition given above based on the results
from our simulations. We note that the df for the group Lasso of Yuan
and Lin (\citeyear{YL06}) requires an initial (least squares) estimator, which is
not available when $p > n$. Thus, their df is not applicable to our problem.

In our simulation example, we compare the adaptive group Lasso with the
group Lasso and ordinary Lasso. Here, the ordinary Lasso estimator is
defined as the
value that minimizes
\[
\|{\mathbf Y}-{\mathbf Z}\bolds\beta_n\|_2^2 + \lambda_n \sum
_{j=1}^p\sum
_{k=1}^{m_n}|\beta_{jk}|.
\]
This simple application of the Lasso does not take into account the
grouping structure in the spline expansions of the components.
The group Lasso
and the adaptive group Lasso estimates are computed using the
algorithm proposed by Yuan and Lin (\citeyear{YL06}). The ordinary Lasso estimates
are computed using the Lars algorithms [Efron et al. (\citeyear{Efronetal04})]. The group
Lasso is used as the initial estimate for the adaptive group Lasso.

We also compare the results from the nonparametric additive modeling
with those from the standard linear regression model with Lasso. We
note that this is not a fair comparison because the generating model is
highly nonlinear. Our purpose is to illustrate that it is necessary to
use nonparametric models when the underlying model deviates
substantially from linear models in the context of variable selection
with high-dimensional data and that model misspecification can lead to
bad selection results.
\begin{Example}\label{example1}
We generate data from the model
\[
y_{i}=f(x_i)+\varepsilon_i \equiv\sum
_{j=1}^{p}f_{j}(x_{ij})+\varepsilon
_{i},\qquad i=1,\ldots, n,
\]
where
\(
f_{1}(t)=5t, f_{2}(t)=3(2t-1)^2, f _{3}(t)=4{\sin(2\pi t)}/{(2-\sin
(2\pi t) )},
\)
\(
f_{4}(t)=6(0.1\sin(2\pi t)+0.2\cos(2\pi t)+
0.3\sin(2\pi t)^2+0.4\cos(2 \pi t)^3+0.5\sin(2\pi t)^3),
\)
and $f_5(t)=\cdots= f_p(t)=0$. Thus, the number of nonzero functions
is $q=4$. This generating model is the same as Example 1 of Lin and
Zhang (\citeyear{LinZhang06}). However, here we use this model in high-dimensional settings.
We consider the cases where $p=1000$ and
three different sample sizes: $n=50, 100$ and $200$.
We use the cubic B-spline with six evenly
distributed knots for all the functions $f_{k}$. The number of
replications in all the simulations is 400.

The covariates are simulated as follows. First, we
generate $w_{i1},\ldots, w_{ip}, u_{i}$, $ u^{\prime}_{i}, v_{i}$
independently from $N(0,1)$ truncated to the interval $[0,1]$,
$i=1,\ldots,n$. Then we set $x_{ik}=(w_{ik}+tu_{i})/(1+t)$
for $k=1,\ldots, 4$ and
$x_{ik}=(w_{ik}+tv_{i})/(1+t)$ for $k=5,\ldots, p$, where the
parameter $t$ controls the amount of correlation among predictors. We
have $\operatorname{Corr}(x_{ik},x_{ij})=t^2/(1+t^2)$, $1\leq j \leq4$, $1\leq k
\leq4$, and $\operatorname{Corr}(x_{ik},x_{ij})=t^2/(1+t^2)$, $4\leq j \leq p$,
$4\leq k \leq p$, but the covariates of the nonzero components and
zero components are independent. We consider $t=0, 1$ in our
simulation. The signal to noise ratio is defined to be
$sd(f)/sd(\epsilon)$. The error term is chosen to be
$\epsilon_{i}\sim N(0, 1.27^2)$ to give a signal-to-noise ratio (SNR)
$3.11:1$. This value is the same as the estimated SNR in the real data
example below, which is the square root of the ratio of the sum of
estimated components squared divided by the sum of residual squared.

The results of 400 Monte Carlo replications are summarized in Table
\ref{table1}. The columns are the mean number of variables selected
(NV), model error (ER), the percentage of replications in which all the
correct additive components are included in the selected model (IN),
and the percentage of replications in which precisely the correct
components are selected (CS). The corresponding standard errors are in
parentheses. The model error is computed as the average of
$n^{-1}\sum_{i=1}^n[\hat{f}(x_i)-f(x_i)]^2$ over the 400 Monte Carlo
replications, where $f$ is the true conditional mean function.

%
%
\begin{sidewaystable}
\tablewidth=\textheight
\tablewidth=\textwidth
\tabcolsep=0pt
\caption{Example 1. Simulation results for the adaptive group Lasso,
group Lasso, ordinary Lasso, and linear model with Lasso, $n=50, 100$
or $200$, $p=1000$.
NV, average number of the variables being selected;
ME, model error; IN, percentage of occasions on which the correct
components are included in the selected model; CS, percentage of
occasions on which correct
components are selected, averaged over 400 replications. Enclosed in
parentheses are the corresponding standard errors. Top panel,
independent predictors; bottom panel, correlated predictors}
\label{table1}
{\fontsize{8.4}{9.8}\selectfont{
\begin{tabular*}{\tablewidth}{@{\extracolsep{4in minus
4in}}lccd{2.3}d{2.3}d{2.3}cd{2.3}
d{2.3}d{2.3}d{2.3}d{2.3}d{2.3}d{2.3}cd{2.3}d{2.3}c@{}}
\hline
& &\multicolumn{4}{c}{\textbf{Adaptive group Lasso}}
& \multicolumn{4}{c}{\textbf{Group Lasso}}
& \multicolumn{4}{c}{\textbf{Ordinary Lasso}}
& \multicolumn{4}{c@{}}{\textbf{Linear mode with Lasso}}\\[-4pt]
& &\multicolumn{4}{c}{\hrulefill}
& \multicolumn{4}{c}{\hrulefill}
& \multicolumn{4}{c}{\hrulefill}
& \multicolumn{4}{r@{}}{\hrulefill}\\
& & \multicolumn{1}{c}{\textbf{NV}} & \multicolumn{1}{c}{\textbf
{ME}} &
\multicolumn{1}{c}{\textbf{IN}} & \multicolumn{1}{c}{\textbf{CS}} &
\multicolumn{1}{c}{\textbf{NV}} &
\multicolumn{1}{c}{\textbf{ME}} & \multicolumn{1}{c}{\textbf{IN}} &
\multicolumn{1}{c}{\textbf{CS}} &
\multicolumn{1}{c}{\textbf{NV}} & \multicolumn{1}{c}{\textbf{ME}} &
\multicolumn{1}{c}{\textbf{IN}} &
\multicolumn{1}{c}{\textbf{CS}} & \multicolumn{1}{c}{\textbf{NV}} &
\multicolumn{1}{c}{\textbf{ME}} &
\multicolumn{1}{c}{\textbf{IN}} & \multicolumn{1}{c@{}}{\textbf
{CS}}\\
\hline
&&\multicolumn{16}{c@{}}{Independent predictors}\\
[2pt]
$n=200$&BIC&4.15&26.72&90.00&80.00 &4.20&27.54&90.00&58.25
&9.73&28.44&95.00&18.00 &3.35&31.89&0.00&0.00 \\
&&(0.43)&(4.13)&(0.30)&(0.41) &(0.43)&(4.45)&(0.30)&(0.54)
&(6.72)&(5.55)&(0.22)&(0.40) &(1.75)&(5.65)&(0.00)&(0.00) \\
&EBIC&4.09&26.64&92.00&81.75 &4.18&27.40&92.00&60.00
&9.58&28.15&95.00&32.50 &3.30&32.08&0.00&0.00\\
&&(0.38)&(4.06)&(0.24)&(0.39) &(0.40)&(4.33)&(0.24)&(0.50)
&(6.81)&(5.25)&(0.22)&(0.47) &(1.86)&(5.69)&(0.00)&(0.00)\\
$n=100$&BIC&4.73&28.26&85.00&70.00 &5.03&29.07&85.00&35.00
&17.25&29.50&82.50&12.00 &6.35&31.57&5.00&0.00\\
&&(1.18)&(5.71)&(0.36)&(0.46) &(1.22)&(6.01)&(0.36)&(0.48)
&(8.72)&(5.89)&(0.38)&(0.44)
&(2.91)&(7.22)&(0.22)&(0.00)\\
&EBIC&4.62&28.07&84.25&74.00 &4.90&28.87&84.25&38.00
&15.93&29.35&84.00&27.75 &5.90&31.53&5.00&0.00\\
&&(0.89)&(5.02)&(0.36)&(0.42) &(1.20)&(5.72)&(0.36)&(0.50)
&(9.06)&(5.25)&(0.36)&(0.45)
&(2.97)&(6.40)&(0.22)&(0.00)\\
$n=50$&BIC&4.75&28.86&80.00&65.00 &5.12&29.97&80.00&32.00
&18.53&30.05&75.00&11.00 &12.53&32.52&22.50&0.00\\
&&(1.22)&(5.72)&(0.41)&(0.48) &(1.29)&(6.15)&(0.41)&(0.48)
&(12.67)&(6.26)&(0.41)&(0.31)
&(3.80)&(8.37)&(0.43)&(0.00)\\
&EBIC&4.69&28.94&78.00&65.00 &5.01&29.82&78.00&36.00
&17.27&30.50&77.50&26.00 &10.33&31.64&20.00&0.00\\
&&(1.98)&(6.48)&(0.40)&(0.48) &(1.21)&(6.11)&(0.40)&(0.49)
&(15.32)&(7.89)&(0.39)&(0.44)
&(3.19)&(8.17)&(0.41)&(0.00)\\
[2pt]
&&\multicolumn{16}{c@{}}{Correlated predictors}\\
[2pt]
$n=200$&BIC&3.20&27.76&66.00&60.00 &3.85&28.12&66.00&30.00
&9.13&28.80&56.00&11.00 &1.08&32.18&0.00&0.00\\
&&(1.27)&(4.74)&(0.46)&(0.50) &(1.49)&(4.76)&(0.46)&(0.46)
&(7.02)&(5.36)&(0.51)&(0.31)
&(0.33)&(8.99)&(0.00)&(0.00)\\
&EBIC&3.23&27.60&68.00&63.00 &3.92&27.85&68.00&31.00
&9.24&28.22&58.00&13.75 &1.30&32.00&0.00&0.00\\
&&(1.24)&(4.34)&(0.45)&(0.49) &(1.68)&(4.50)&(0.45)&(0.48)
&(7.18)&(5.30)&(0.52)&(0.44)
&(1.60)&(8.92)&(0.00)&(0.00)\\
$n=100$&BIC&2.88&27.88&60.00&56.00 &3.28&28.33&60.00&22.00
&8.80&28.97&52.00&8.00 &1.00&32.24&0.00&0.00\\
&&(1.91)&(4.88)&(0.50)&(0.56) &(1.96)&(4.92)&(0.50)&(0.42)
&(10.22)&(5.45)&(0.44)&(0.26)
&(0.00)&(9.20)&(0.00)&(0.00)\\
&EBIC&3.04&27.78&61.75&58.00 &3.44&28.16&61.75&24.00
&9.06&28.55&54.00&10.00 &1.00&32.09&0.00&0.00\\
&&(1.46)&(4.85)&(0.49)&(0.54) &(1.52)&(4.90)&(0.49)&(0.43)
&(11.24)&(5.42)&(0.46)&(0.28)
&(0.00)&(8.98)&(0.00)&(0.00)\\
$n=50$&BIC&2.50&28.36&48.50&38.00 &3.10&29.37&48.50&20.00
&8.01&30.48&30.00&5.00 &1.00&33.28&0.00&0.00\\
&&(1.64)&(5.32)&(0.50)&(0.55) &(1.78)&(5.98)&(0.50)&(0.41)
&(11.42)&(6.77)&(0.46)&(0.23)
&(0.00)&(9.42)&(0.00)&(0.00)\\
&EBIC&2.48&28.57&48.00&38.00 &3.07&30.13&48.00&18.00
&8.24&30.89&32.00&6.00 &1.00&33.25&0.00&0.00\\
&&(1.62)&(5.51)&(0.51)&(0.55) &(1.76)&(7.60)&(0.51)&(0.40)
&(11.46)&(6.40)&(0.48)&(0.24)
&(0.00)&(9.38)&(0.00)&(0.00)\\
\hline
\end{tabular*}
}}
\end{sidewaystable}

Table \ref{table1} shows that the adaptive group Lasso selects all the
nonzero components (IN) and selects exactly the correct model (CS) more
frequently than the other methods do. For example, with the BIC and
$n=200$, the percentage of correct selections (CS) by the adaptive
group Lasso ranges from 65.25\% to 81\%, which is much higher than the
ranges 30--57.75\% for the group Lasso and 12--15.75\% for the ordinary
Lasso. The adaptive group Lasso and group Lasso perform better than the
ordinary Lasso in all of the experiments, which illustrates the
importance of taking account of the group structure of the coefficients
of the spline expansion. Correlation among covariates increases the
difficulty of component selection, so it is not surprising that all
methods perform better with independent covariates than with correlated
ones. The percentage of correct selections increases as the sample size
increases. The linear model with Lasso never selects the correct model.
This illustrates the poor results that can be produced by a linear
model when the true conditional mean function is nonlinear.

Table \ref{table1} also shows that the model error (ME) of the group
Lasso is only slightly larger than that of the adaptive group Lasso.
The models selected by the group Lasso nest and, therefore, have more
estimated coefficients than the models selected by the adaptive group
Lasso. Therefore, the group Lasso estimators of the conditional mean
function have a larger variance and larger ME. The differences between
the MEs of the two methods are small, however, because as can be seen
from the NV column, the models selected by the group Lasso in our
experiments have only slightly more estimated coefficients than the
models selected by the adaptive group Lasso.
\end{Example}
\begin{Example}\label{example2}
We now compare the adaptive group Lasso with the COSSO [Lin and Zhang
(\citeyear{LinZhang06})]. This comparison is suggested to us by the Associate Editor.
Because the COSSO algorithm only works for the case when $p$ is smaller
than $n$, we use the same set-up as in Example 1 of Lin and Zhang
(\citeyear{LinZhang06}). In this example, the generating model is as in (\ref{GModela})
with 4 nonzero components.
Let $X_{j}=(W_{j}+tU)/(1+t)$, $j=1,\ldots,p$, where $W_{1},\ldots,
W_{p}$ and $U$ are i.i.d. from $N(0,1)$, truncated to the interval
$[0,1]$. Therefore, corr$(X_{j}, X_{k})=t^2/(1+t^2)$ for $j \ne k$.
The random error term $\epsilon\sim N(0, 1.32^2)$. The SNR is 3:1.
We consider three different sample sizes $n=50, 100$ or $200$ and
three different number of predictors $p=10, 20$ or $50$.
The COSSO estimator is computed using the Matlab software which is
publicly available at
\url{http://www4.stat.ncsu.edu/\textasciitilde hzhang/cosso.html}.

%
%
\begin{sidewaystable}
\tablewidth=\textheight
\tablewidth=\textwidth
\caption{Example 2. Simulation results comparing the adaptive group
Lasso and COSSO. $n=50, 100$ or $200$, $p=10, 20$ or $50$.
NV, average number of the variables being selected;
ME, model error; IN, percentage of occasions on which all the correct
components are included in the selected model; CS, percentage of
occasions on which correct
components are selected, averaged over 400 replications. Enclosed in
parentheses are the corresponding standard errors}
\label{table2}
\begin{tabular*}{\tablewidth}{@{\extracolsep{\fill
}}lcccd{3.3}d{2.3}ccd{3.3}d{2.3}ccd{2.3}d{2.3}@{}}
\hline
&&\multicolumn{4}{c}{$\bolds{p=10}$} &
\multicolumn{4}{c}{$\bolds{p=20}$} &
\multicolumn{4}{c@{}}{$\bolds{p=50}$}\\[-4pt]
&&\multicolumn{4}{c}{\hrulefill} &
\multicolumn{4}{c}{\hrulefill} &
\multicolumn{4}{c@{}}{\hrulefill}\\
&&\multicolumn{1}{c}{\textbf{NV}} & \multicolumn{1}{c}{\textbf{ME}}
& \multicolumn{1}{c}{\textbf{IN}} & \multicolumn{1}{c}{\textbf{CS}}
& \multicolumn{1}{c}{\textbf{NV}} & \multicolumn{1}{c}{\textbf{ME}}
& \multicolumn{1}{c}{\textbf{IN}} & \multicolumn{1}{c}{\textbf{CS}}
& \multicolumn{1}{c}{\textbf{NV}} & \multicolumn{1}{c}{\textbf{ME}}
& \multicolumn{1}{c}{\textbf{IN}} & \multicolumn{1}{c@{}}{\textbf{CS}}\\
\hline
&&\multicolumn{12}{c@{}}{Independent predictors}\\
[2pt]
$n=200$& AGLasso(BIC)&4.02&0.27&100.00&98.00 &4.01&0.34&96.00&92.00
&4.10&0.88&98.00&90.00\\
&&(0.14)&(0.10)&(0.00)&(0.14) &(0.40)&(0.10)&(0.20)&(0.27)
&(0.39)&(0.19)&(0.14)&(0.30)\\
&AGLasso(EBIC)&4.02&0.27&100.00&99.00 &4.05&0.32&100.00&94.00
&4.08&0.87&98.00&90.00\\
&&(0.14)&(0.09)&(0.00)&(0.10) &(0.22)&(0.09)&(0.00)&(0.24)
&(0.30)&(0.16)&(0.14)&(0.30)\\
& COSSO(5CV)&4.06&0.29&100.00&98.00 &4.10&0.37&100.00&92.00
&4.49&1.53&94.00&84.00 \\
&&(0.24)&(0.07)&(0.00)&(0.14) &(0.39)&(0.11)&(0.00)&(0.27)
&(1.10)&(0.86)&(0.24)&(0.37)\\
[2pt]
$n=100$& AGLasso(BIC)&4.06&0.56&99.00&90.00 &4.11&0.63&98.00&87.00
&4.27&1.04&93.00&81.00\\
&&(0.24)&(0.19)&(0.10)&(0.30) &(0.42)&(0.26)&(0.14)&(0.34)
&(0.58)&(0.64)&(0.26)&(0.39)\\
& AGLasso(EBIC)&4.06&0.54&99.00&91.00 &4.10&0.59&98.00&89.00
&4.22&1.01&93.00&83.00\\
&&(0.24)&(0.21)&(0.10)&(0.31) &(0.39)&(0.22)&(0.14)&(0.31)
&(0.56)&(0.60)&(0.26)&(0.38)\\
& COSSO(5CV)&4.17&0.53&96.00&89.00 &4.18&1.04&83.00&63.00
&4.89&6.63&30.00&11.00\\
&&(0.62)&(0.19)&(0.20)&(0.31) &(0.96)&(0.64)&(0.38)&(0.49)
&(1.50)&(1.29)&(0.46)&(0.31)\\
[2pt]
$n=50$&AGLasso(BIC)&4.18&0.72&98.00&84.00 &4.25&0.99&96.00&79.00
&4.30&1.06&90.00&71.00\\
&&(0.66)&(0.56)&(0.14)&(0.36) &(0.72)&(0.60)&(0.20)&(0.41)
&(0.89)&(0.68)&(0.30)&(0.46)\\
&AGLasso(EBIC)&4.16&0.70&98.00&84.00 &4.24&1.02&94.00&78.00
&4.27&1.04&92.00&73.00\\
&&(0.64)&(0.52)&(0.14)&(0.36) &(0.70)&(0.62)&(0.20)&(0.42)
&(0.86)&(0.64)&(0.27)&(0.45)\\
&COSSO(5CV)&4.41&1.77&61.00&58.00 &5.06&5.53&33.00&20.00
&5.96&7.60&8.00&0.00\\
&&(1.08)&(1.35)&(0.46)&(0.42) &(1.54)&(1.88)&(0.47)&(0.40)
&(2.20)&(2.07)&(0.27)&(0.00)\\
\hline
\end{tabular*}
\end{sidewaystable}

\setcounter{table}{1}
%
%
\begin{sidewaystable}
\tablewidth=\textheight
\tablewidth=\textwidth
\caption{(Continued)}
\begin{tabular*}{\tablewidth}{@{\extracolsep{\fill
}}lcccd{3.3}d{2.3}ccd{3.3}d{2.3}ccd{2.3}d{2.3}@{}}
\hline
&&\multicolumn{4}{c}{$\bolds{p=10}$} &
\multicolumn{4}{c}{$\bolds{p=20}$} &
\multicolumn{4}{c@{}}{$\bolds{p=50}$}\\[-4pt]
&&\multicolumn{4}{c}{\hrulefill} &
\multicolumn{4}{c}{\hrulefill} &
\multicolumn{4}{c@{}}{\hrulefill}\\
&&\multicolumn{1}{c}{\textbf{NV}} & \multicolumn{1}{c}{\textbf{ME}}
& \multicolumn{1}{c}{\textbf{IN}} & \multicolumn{1}{c}{\textbf{CS}}
& \multicolumn{1}{c}{\textbf{NV}} & \multicolumn{1}{c}{\textbf{ME}}
& \multicolumn{1}{c}{\textbf{IN}} & \multicolumn{1}{c}{\textbf{CS}}
& \multicolumn{1}{c}{\textbf{NV}} & \multicolumn{1}{c}{\textbf{ME}}
& \multicolumn{1}{c}{\textbf{IN}} & \multicolumn{1}{c@{}}{\textbf{CS}}\\
\hline
&&\multicolumn{12}{c@{}}{Correlated predictors}\\
[2pt]
$n=200$&AGLasso(BIC)&3.75&0.49&82.00&70.00 &3.71&1.20&75.00&66.00
&3.50&1.68&68.00&62.00\\
&&(0.61)&(0.14)&(0.39)&(0.46) &(0.68)&(0.89)&(0.41)&(0.46)
&(0.92)&(1.29)&(0.45)&(0.49)\\
&AGLasso(EBIC)&3.75&0.49&82.00&70.00 &3.73&1.18&75.00&68.00
&3.58&1.60&70.00&65.00\\
&&(0.61)&(0.14)&(0.39)&(0.46) &(0.65)&(0.88)&(0.41)&(0.45)
&(0.84)&(1.27)&(0.46)&(0.46)\\
&COSSO(5CV)&3.70&0.53&69.00&41.00 &3.89&1.24&57.00&36.00
&4.11&1.76&41.00&16.00\\
&&(0.58)&(0.17)&(0.46)&(0.49) &(0.60)&(0.90)&(0.50)&(0.48)
&(0.86)&(1.33)&(0.49)&(0.37)\\
[2pt]
$n=100$&AGLasso(BIC)&3.72&1.40&78.00&68.00 &3.68&1.78&70.00&64.00
&3.02&3.07&63.00&59.00\\
&&(0.66)&(0.70)&(0.40)&(0.45) &(0.74)&(1.15)&(0.46)&(0.48)
&(1.58)&(2.37)&(0.49)&(0.51)\\
&AGLasso(EBIC)&3.70&1.46&75.00&66.00 &3.71&1.74&72.00&64.00
&3.20&2.98&65.00&60.00\\
&&(0.72)&(0.78)&(0.41)&(0.46) &(0.68)&(1.06)&(0.42)&(0.48)
&(1.42)&(1.96)&(0.46)&(0.50)\\
&COSSO(5CV)&3.98&1.42&41.00&26.00 &4.14&1.76&30.00&6.00
&4.24&6.88&8.00&0.00\\
&&(0.64)&(0.74)&(0.49)&(0.42) &(2.27)&(1.11)&(0.46)&(0.24)
&(2.96)&(2.91)&(0.27)&(0.00)\\
[2pt]
$n=50$&AGLasso(BIC)&3.30&2.26&70.00&62.00 &3.06&3.02&65.00&60.00
&2.87&4.01&52.00&42.00\\
&&(1.16)&(1.09)&(0.46)&(0.49) &(1.52)&(2.14)&(0.46)&(0.50)
&(1.56)&(3.69)&(0.44)&(0.52)\\
&AGLasso(EBIC)&3.32&2.20&70.00&64.00 &3.10&3.01&68.00&62.00
&2.90&3.88&50.00&42.00\\
&&(1.14)&(1.06)&(0.46)&(0.48) &(1.51)&(2.12)&(0.45)&(0.49)
&(1.54)&(3.62)&(0.42)&(0.52)\\
&COSSO(5CV)&4.14&3.77&25.00&6.00 &4.20&6.98&5.00&0.00
&4.90&9.93&1.00&0.00\\
&&(2.25)&(2.02)&(0.44)&(0.24) &(2.88)&(2.82)&(0.22)&(0.00)
&(3.30)&(4.08)&(0.10)&(0.00)\\
\hline
\end{tabular*}
\end{sidewaystable}

The COSSO procedure uses either generalized cross-validation or
5-fold cross-validation. Based the simulation results of Lin and Zhang
(\citeyear{LinZhang06}) and our own simulations, the COSSO with 5-fold cross-validation
has better selection performance. Thus, we compare the adaptive group Lasso
with BIC or EBIC with the COSSO with 5-fold cross-validation. The
results are given in Table \ref{table2}. For independent predictors,
when $n=200$ and $p=10, 20$ or $50$, the adaptive group Lasso and
COSSO have similar performance in terms of selection accuracy and model
error. However, for smaller $n$ and
larger $p$, the adaptive group Lasso does significantly better. For
example, for $n=100$ and $p=50$, the percentage of correct selection for
the adaptive group Lasso is 81--83\%, but it is only 11\% for the
COSSO. The model error of the adaptive group Lasso is
similar to or smaller than that of the COSSO. In several experiments,
the model error of the COSSO is 2 to more than 7 times larger than that
of the adaptive group Lasso. It is interesting to note that when $n=50$
and $p=20$ or $50$, the adaptive group Lasso still does a descent job
in selecting the correct model, but the COSSO does poorly in these two
cases. In particular, for $n=50$ and $p=50$, the COSSO did not select
the exact correct model in all the simulation runs. For dependent
predictors, the comparison is even mode favorable to the adaptive group
Lasso, which performs significantly better than COSSO in terms of both
model error and selection accuracy in all the cases.
\end{Example}

\section{Data example}\label{sec5}

We use the data set reported in Scheetz et al. (\citeyear{Scheetzetal06}) to
illustrate the application of the proposed method in
high-dimensional settings. For this data set,
120 twelve-week old male rats were selected
for tissue harvesting from the eyes and for microarray analysis. The
microarrays used to analyze the RNA from the eyes of these 
animals contain over 31,042 different probe sets (Affymetric
GeneChip Rat Genome 230 2.0 Array). The intensity values were
normalized using the robust multi-chip averaging method [Irizzary et
al. (\citeyear{Irizarryetal03})] method to obtain summary expression values
for each probe set. Gene expression levels were analyzed on a
logarithmic scale.

We are interested in finding the genes that are related to the gene
TRIM32. This gene was recently found to cause Bardet--Biedl syndrome
[Chiang et al. (\citeyear{Chiangetal06})], which is a genetically heterogeneous disease of
multiple organ systems including the retina. Although over 30,000 probe
sets are represented on the Rat Genome 230 2.0 Array, many of them are
not expressed in the eye tissue and initial screening using correlation
shows that most probe sets have very low correlation with TRIM32. In
addition, we are expecting only a small number of genes to be related
to TRIM32. Therefore, we use 500 probe sets that are expressed in the
eye and have highest marginal correlation in the analysis. Thus, the
sample size is $n=120$ (i.e., there are 120 arrays from 120 rats) and $p=500$.
It is expected that only a few genes are related to TRIM32.
Therefore, this is a sparse, high-dimensional regression problem.

We use the nonparametric additive model to model the relation between
the expression of TRIM32 and those of the 500 genes.
We estimate model (\ref{GamA})
using the ordinary Lasso, group Lasso, and adaptive
group Lasso for the nonparametric additive model. To compare the results
of the nonparametric additive model with that of the linear regression
model, we also analyzed the data using the linear regression model with
Lasso. We scale the covariates so that their values are between 0 and 1
and use cubic splines with six evenly distributed knots to estimate the
additive components. The penalty parameters in all the methods are
chosen using the BIC or EBIC as in the simulation study. Table \ref
{table3} lists the probes selected by
the group Lasso and the adaptive group Lasso, indicated by the check
signs. Table \ref{table4}
shows the number of variables, the residual sums of squares obtained
with each estimation method. For the ordinary Lasso with the spline
expansion, a variable is
considered to be selected if any of the estimated coefficients of
the spline approximation to its additive component are nonzero.
Depending on whether BIC or EBIC is used, the
group Lasso selects 16--17 variables, the adaptive group Lasso selects
15 variables and the ordinary Lasso with the spline expansion selects
94--97 variables, the linear model selects 8--14 variables. Table \ref
{table4} shows that the adaptive group Lasso does better than the other
methods in terms of residual sum of squares (RSS). We have also
examined the plots (not shown) of the
estimated additive components obtained with the group Lasso and the
adaptive group Lasso, respectively. Most are highly nonlinear,
confirming the need for taking into account nonlinearity.

%
%
\begin{table}
\tabcolsep=0pt
\caption{Probe sets selected by the group Lasso and
the adaptive group Lasso in the data example using BIC or EBIC for
penalty parameter selection. GL, group Lasso; AGL, adaptive group
Lasso; Linear, linear model with Lasso}\label{table3}
\begin{tabular*}{\tablewidth}{@{\extracolsep{\fill}}lcccccc@{}}
\hline
\textbf{Probes} & \textbf{GL(BIC)} & \textbf{AGL(BIC)}
& \textbf{Linear(BIC)} & \textbf{GL(EBIC)} & \textbf{AGL(EBIC)}
& \textbf{Linear(EBIC)}\\
\hline
$1389584\_at$ & $\surd$ & $\surd$&$\surd$& $\surd$ &$\surd$
&$\surd$ \\
$1383673\_at$ & $\surd$ & $\surd$& $\surd$&$\surd$& $\surd$&$\surd
$\\
$1379971\_at$ & $\surd$ & $\surd$& $\surd$&$\surd$ & $\surd
$&$\surd$\\
$1374106\_at$ & $\surd$ & &$\surd$&$\surd$ & &$\surd$\\
$1393817\_at$ & $\surd$ & $\surd$& $\surd$&$\surd$&$\surd$&\\
$1373776\_at$ & $\surd$ & $\surd$& $\surd$&$\surd$&$\surd$&\\
$1377187\_at$ & $\surd$ & $\surd$& $\surd$&$\surd$&$\surd$&\\
$1393955\_at$ & $\surd$ & $\surd$&$\surd$ &$\surd$&$\surd$&\\
$1393684\_at$ & $\surd$ & $\surd$& &$\surd$&$\surd$&\\
$1381515\_at$ & $\surd$ & $\surd$& &$\surd$&$\surd$&\\
$1382835\_at$ & $\surd$ & $\surd$&$\surd$ &$\surd$&$\surd$&\\
$1385944\_at$ & $\surd$ & $\surd$&$\surd$ &$\surd$&$\surd$&\\
$1382263\_at$ & $\surd$ & $\surd$& $\surd$&$\surd$&$\surd$&$\surd
$\\
$1380033\_at$ & $\surd$ & $\surd$& &$\surd$&$\surd$&\\
$1398594\_at$ & $\surd$ & &$\surd$&&&$\surd$\\
$1376744\_at$ & $\surd$ & $\surd$& &$\surd$&$\surd$&\\
$1382633\_at$ & $\surd$ & $\surd$& &$\surd$&$\surd$&\\
$1383110\_at$ & & & $\surd$&&&$\surd$\\
$1386683\_at$&&&$\surd$&&&$\surd$\\
\hline
\end{tabular*}
\end{table}

%
%
\begin{table}[b]
\caption{Analysis results for the data example. No. of probes, the
number of probe sets selected; RSS,
the residual sum of squares of the fitted model}\label{table4}
\begin{tabular*}{\tablewidth}{@{\extracolsep{\fill}}l@{}cccc@{}}
\hline
& \multicolumn{2}{c}{\textbf{BIC}} &
\multicolumn{2}{c@{}}{\textbf{EBIC}}\\[-4pt]
& \multicolumn{2}{c}{\hrulefill} & \multicolumn{2}{c@{}}{\hrulefill
}\\
& \textbf{No. of probe sets} & \textbf{RSS} & \textbf{No. of probe sets}
& \textbf{RSS}\\
\hline
Adaptive group Lasso& 15 & 1.52e--03&15&1.52e--03\\
Group Lasso& 17 & 3.24e--03&16&3.40e--03\\
Ordinary Lasso & 97 & 2.96e--07&94&8.10e--08\\
Linear regression with Lasso & 14 & 2.62e--03&\phantom{0}8&3.75e--03\\
\hline
\end{tabular*}
\end{table}

In order to evaluate the performance of the methods, we
use cross-validation and compare the prediction mean square errors
(PEs). We randomly partition the data into 6
subsets, each set consisting of 20 observations. We then fit the
model with 5 subsets as training set and calculate the PE
for the remaining set which we consider as test set. We repeat this
process 6 times, considering one of the 6 subsets as test set
every time. We compute the average of the numbers of probes selected and
the prediction errors of these 6 calculations. Then we replicate this
process 400 times (this is suggested to us by the Associate Editor).
Table \ref{table5} gives the average values over 400 replications.
The adaptive group Lasso has smaller average prediction error than the
group Lasso, the ordinary Lasso and the linear regression with Lasso.
The ordinary Lasso selects far more probe sets than the other
approaches, but this does not lead to better prediction performance.
Therefore, in this example, the adaptive group Lasso provides the
investigator a more targeted list of probe sets, which can serve as a
starting point for further study.

%
%
\begin{table}
\tabcolsep=0pt
\caption{Comparison of adaptive group Lasso, group Lasso, ordinary
Lasso, and linear regression model with Lasso for the data example.
ANP, the average number of probe sets selected averaged across 400
replications; PE, the average of prediction mean square errors for
the test set}\label{table5}
\begin{tabular*}{\tablewidth}{@{\extracolsep{4in minus
4in}}ld{2.3}cd{2.3}cd{2.3}cd{2.3}c@{}}
\hline
&\multicolumn{2}{c}{\textbf{Adaptive}} &
&
&&&
\multicolumn{2}{c}{\textbf{Linear}}\\
&\multicolumn{2}{c}{\textbf{group Lasso}} &
\multicolumn{2}{c}{\textbf{Group Lasso}} &
\multicolumn{2}{c}{\textbf{Ordinary Lasso}} &
\multicolumn{2}{c}{\textbf{model with Lasso}}\\[-4pt]
&\multicolumn{2}{c}{\hrulefill} &
\multicolumn{2}{c}{\hrulefill} &
\multicolumn{2}{c}{\hrulefill} &
\multicolumn{2}{c@{}}{\hrulefill}\\
& \multicolumn{1}{c}{\textbf{ANP}} & \multicolumn{1}{c}{\textbf{PE}}
& \multicolumn{1}{c}{\textbf{ANP}} & \multicolumn{1}{c}{\textbf{PE}}
& \multicolumn{1}{c}{\textbf{ANP}} & \multicolumn{1}{c}{\textbf{PE}}
& \multicolumn{1}{c}{\textbf{ANP}} & \multicolumn{1}{c@{}}{\textbf
{PE}}\\
\hline
BIC&15.75&1.86e--02&16.45&2.89e--02&78.48&1.40e--02&9.25&2.26e--02\\
&(0.85)&(0.47e--02)&(0.88)&(0.49e--02)&(3.62)&(0.90e--02)&(0.88)&(1.41e--2)\phantom
{0}\\
EBIC&15.55&1.78e--02&16.75&1.99e--02&80.00&1.23e--02&9.15&2.03e--02\\
&(0.82)&(0.42e--02)&(0.84)&(0.47e--02)&(3.50)&(0.89e--02)&(0.86)&(1.39e--02)\\
\hline
\end{tabular*}
\end{table}

It is of interest to compare the selection results from the adaptive
group Lasso and the linear regression model with Lasso.
The adaptive group Lasso and the linear model with Lasso select
different sets of genes. When the penalty parameter is chosen with the
BIC, the adaptive group Lasso selects 5 genes that are not selected by
the linear model with Lasso. In addition, the linear model with Lasso
selects 5 genes that are not selected by the adaptive group Lasso. When
the penalty parameter is selected with the EBIC, the adaptive group
Lasso selects 10 genes that are not selected by the linear model with
Lasso. The estimated effects of many of the genes are nonlinear, and
the Monte Carlo results of Section \ref{sec4} show that the
performance of the
linear model with Lasso can be very poor in the presence of
nonlinearity. Therefore, we interpret the differences between the gene
selections of the adaptive group Lasso and the linear model with Lasso
as evidence that the selections produced by the linear model are misleading.

\section{Concluding remarks}\label{sec6}

In this paper, we propose to use the adaptive group Lasso
for variable selection in nonparametric additive models in sparse,
high-dimensional settings. A key requirement for the adaptive group
Lasso to be selection consistent is that the initial estimator is
estimation consistent and selects all the important components with
high probability. In low-dimensional settings, finding an initial
consistent estimator is relatively easy and can be achieved by many
well-established approaches such as the additive spline estimators.
However, in high-dimensional settings, finding an initial consistent
estimator is difficult. Under the conditions stated in Theorem
\ref{ThmA}, the group Lasso
is shown to be consistent and selects all the important components.
Thus the group Lasso can be used as the initial estimator in the adaptive
Lasso to achieve selection consistency.
Following model selection, oracle-efficient, asymptotically normal
estimators of the nonzero components can be obtained by using existing
methods. Our simulation results indicate that our procedure works well
for variable selection in the models considered. Therefore, the
adaptive group Lasso is a useful approach for variable selection and
estimation in sparse, high-dimensional nonparametric additive models.

Our theoretical results are concerned with a fixed sequence of penalty
parameters, which are not applicable to the case where the penalty
parameters are selected based on data driven procedures such as the
BIC. This is an important and challenging problem that deserves further
investigation, but is beyond the scope of this paper.
We have only considered linear nonparametric additive models. The
adaptive group Lasso can be applied to generalized nonparametric
additive models, such as the generalized logistic nonparametric
additive model and other nonparametric models with high-dimensional
data. However,
more work is needed to understand the properties of this approach in
those more complicated models.

\begin{appendix}\label{app}
\section*{Appendix: Proofs}

We first prove the following lemmas. Denote the centered versions of
$\mathcal{S}_n$ by
\[
\mathcal{S}_{nj}^0 = \Biggl\{f_{nj}\dvtx f_{nj}(x) =
\sum_{k=1}^{m_n} b_{jk}\psi_{k}(x), (\beta_{j1}, \ldots,
\beta_{jm_n}) \in\mathbb{R}^{m_n} \Biggr\},\qquad 1 \le j \le p,
\]
where $\psi_k$'s are the centered spline bases defined in (\ref{CenterA}).
\begin{lemma}
\label{LemA} Suppose that $f \in\mathcal{F}$ and $\mathrm
{E}f(X_j)=0$. Then under
\textup{(A3)} and \textup{(A4)}, there
exists an $f_n \in\mathcal{S}_{nj}^0$ satisfying
\[
\|f_{n}-f\|_{2} = O_p(m_n^{-d}+m_n^{1/2} n^{-1/2}).
\]
In particular, if we choose $m_n =O(n^{1/(2d+1)})$, then
\[
\|f_{n}-f\|_{2} = O_p(m_n^{-d})=O_p\bigl(n^{-d/(2d+1)}\bigr).
\]
\end{lemma}
\begin{pf}
By (A4), for $f \in\mathcal{F}$, there
is an $f_n^* \in\mathcal{S}_n$ such that $\|f-f_{n}^*\|_{2} =
O(m_n^{-d})$. Let $f_n = f_n^*- n^{-1}\sum_{i=1}^n f_n^*(X_{ij})$.
Then $f_n \in\mathcal{S}^0_{nj}$ and $|f_n-f| \le|f_n^*-f| + |P_n
f_n^*|$, where $P_n$ is the empirical measure of i.i.d. random
variables $X_{1j}, \ldots, X_{nj}$. Consider
\[
P_n f_n^*=(P_n - P)f_n^* + P(f_n^*-f).
\]
Here, we use the linear functional notation, for example, $Pf = \int f dP$,
where $P$ is the probability measure of $X_{1j}$.
For any $\varepsilon> 0$, the bracketing number $N_{[\cdot]}(\varepsilon,
\mathcal{S}_{nj}^0, L_2(P)) $ of $\mathcal{S}_{nj}^0$ satisfies\vspace*{2pt} $\log
N_{[\cdot]}(\varepsilon
, \mathcal{S}_{nj}^0, L_2(P)) \le c_1 m_n \log(1/\varepsilon)$ for some
constant $c_1 > 0$ [Shen and Wong (\citeyear{SW94}), page 597]. Thus, by the maximal
inequality; see, for example, van der Vaart (\citeyear{Vaart1998}, page 288),
$(P_n -P)f_n^* = O_p(n^{-1/2}m_n^{1/2})$. By (A4), $|P(f_n^*-f)|
\le C_2 \|f_n^*-f\|_{2}=O(m_{n}^{-d})$ for some constant $C_2 > 0$. The
lemma follows from the triangle inequality.
\end{pf}
\begin{lemma}
\label{LemC}
Suppose that conditions \textup{(A2)} and \textup{(A4)} hold. Let
\[
T_{jk} = n^{-1/2}m_n^{1/2} \sum_{i=1}^n \psi_k(X_{ij})\varepsilon
_i,\qquad 1 \le j
\le p, 1\le k \le m_n,
\]
and $T_n = {\max_{1\le j\le p,1\le k\le m_n}}|T_{jk}|$. Then
\begin{eqnarray*}
\mathrm{E}(T_n) &\le& C_1 n^{-1/2}m_n^{1/2}
\sqrt{\log(p m_n)} \bigl(\sqrt{ 2 C_2 m_n^{-1} n \log(pm_n)}\\
&&\hspace*{113.2pt}{} + 4
\log(2
p m_n)+ C_2 n m_n^{-1} \bigr)^{1/2},
\end{eqnarray*}
where $C_1$ and $C_2$ are two positive constants. In particular, when $
m_n \log(pm_n)/\break n \rightarrow0$,
\[
\mathrm{E}(T_n) 
=O(1) \sqrt{\log(pm_n)}.
\]
\end{lemma}
\begin{pf}
Let $s_{njk}^2 = \sum_{i=1}^n \psi_k^2(X_{ij})$.\vspace*{-2pt}
Conditional on $X_{ij}$'s, $T_{jk}$'s are sub-Gaussian.
Let $s_n^2 = \max_{1\le j \le p, 1\le k \le m_n} s_{njk}^2$.
By (A2) and the maximal inequality for sub-Gaussian random variables
[van der Vaart and Wellner (\citeyear{VaartWellner96}), Lemmas~2.2.1 and 2.2.2],
\[
\mathrm{E}\Bigl({\max_{1\le j\le p,1\le k\le m_n}}|T_{jk}| \big|\{X_{ij},
1\le
i\le n, 1\le j \le p\} \Bigr) \le C_1 n^{-1/2}m_n^{1/2} s_n
\sqrt{\log(p m_n)}.
\]
Therefore,
%
%
\begin{equation}
\label{BsAa}
\mathrm{E}\Bigl({\max_{1\le j\le p,1\le k\le m_n}}|T_{jk}| \Bigr) \le C_1
n^{-1/2}m_n^{1/2}
\sqrt{\log(p m_n)} \mathrm{E}(s_n),
\end{equation}
where $C_1 > 0$ is a constant. By (A4) and the properties of B-splines,
%
%
\begin{equation}
\label{BsA}
|\psi_k(X_{ij})|\le|\phi_k(X_{ij})| +|\bar{\phi}_{jk}| \le2
\quad\mbox{and}\quad \mathrm{E}(\psi_k(X_{ij}))^2 \le C_2 m_n^{-1}
\end{equation}
for a constant $C_2 > 0$,
for every $ 1\le j \le p$ and $1\le k \le m_n$.
By (\ref{BsA}),
%
%
\begin{equation}
\label{BsB}
\sum_{i=1}^n \mathrm{E}[\psi_k^2(X_{ij})-\mathrm{E}\psi
_k^2(X_{ij})]^2 \le4 C_2
n m_n^{-1}
\end{equation}
and
%
%
\begin{equation}
\label{BsC}
\max_{1\le j\le p, 1\le k \le m_n} \sum_{i=1}^n \mathrm{E}\psi_k^2(X_{ij})
\le C_2 n m_n^{-1}.
\end{equation}
By Lemma A.1 of van de Geer (\citeyear{vandeGeer08}), (\ref{BsA}) and (\ref{BsB}) imply
\begin{eqnarray*}
&&\mathrm{E}\Biggl({\max_{1\le j\le p, 1\le k \le m_n} }\Biggl|\sum_{i=1}^n \{
\psi
_{k}^2(X_{ij})-\mathrm{E}\psi_{k}^2(X_{ij})\} \Biggr| \Biggr)\\
&&\qquad
\le\sqrt{ 2 C_2 m_n^{-1} n \log(pm_n)} + 4 \log(2 p m_n).
\end{eqnarray*}
Therefore, by (\ref{BsC}) and the triangle inequality,
\[
\mathrm{E}s_n^2 \le\sqrt{ 2 C_2 m_n^{-1} n \log(pm_n)} + 4 \log(2
p m_n)
+ C_2 n m_n^{-1}.
\]
Now since $ \mathrm{E}s_n \le( \mathrm{E}s_n^2)^{1/2}$, we have
%
%
\begin{equation}
\label{BsD}
\mathrm{E}s_n \le\bigl(\sqrt{ 2 C_2 m_n^{-1} n \log(pm_n)} + 4 \log(2 p m_n)
+ C_2 n m_n^{-1} \bigr)^{1/2}.
\end{equation}
The lemma follows from (\ref{BsAa}) and (\ref{BsD}).
\end{pf}

Denote
\[
\bolds\beta_{A} = (\bolds\beta_j', j \in A)' \quad\mbox{and}\quad
{\mathbf Z}_A=({\mathbf Z}_j, j \in A).
\]
Here, $\bolds\beta_{A}$ is an $|A|m_n \times1$ vector and ${\mathbf
Z}_A$ is an
$n \times|A|m_n$ matrix.
Let $ {\mathbf C}_A = {\mathbf Z}_A'{\mathbf Z}_A/n. $ When $A=\{1,
\ldots, p\}$, we
simply write ${\mathbf C}={\mathbf Z}'{\mathbf Z}/n$.
Let $\rho_{\min}({\mathbf C}_A)$ and $\rho_{\max}({\mathbf C}_A)$
be the
minimum and maximum eigenvalues of ${\mathbf C}_A$, respectively.
\begin{lemma}
\label{LemPB} Let $m_n = O(n^{\gamma})$ where $0 < \gamma< 0.5$.
Suppose that $|A|$ is bounded by a fixed constant independent of $n$
and $p$. Let $h\equiv h_n \asymp m_n^{-1}$. Then under \textup{(A3)}
and \textup{(A4)},
with probability converging to one,
\[
c_1h_n \le\rho_{\min}({\mathbf C}_A) \le\rho_{\max}({\mathbf
C}_A) \le c_2 h_n,
\]
where $c_1$ and $c_2$ are two positive constants.
\end{lemma}
\begin{pf}
Without loss of generality, suppose
$A=\{1, \ldots, k\}$. Then ${\mathbf Z}_A = ({\mathbf Z}_1$, $\ldots,
{\mathbf Z}_q)$. Let
${\mathbf b}
=({\mathbf b}_1', \ldots, {\mathbf b}_q')'$, where ${\mathbf b}_j \in R^{m_n}$.
By Lemma 3 of Stone (\citeyear{Stone85}),
\[
\|{\mathbf Z}_1 {\mathbf b}_1 + \cdots+ {\mathbf Z}_q {\mathbf b}_q\|
_2 \ge c_3 (\|{\mathbf Z}_1 {\mathbf b}_1\|
_2 +
\cdots+ \|{\mathbf Z}_q{\mathbf b}_q\|_2)
\]
for a certain constant $c_3 > 0$.
By the triangle inequality,
\[
\|{\mathbf Z}_1 {\mathbf b}_1 + \cdots+ {\mathbf Z}_q {\mathbf b}_q\|
_2 \le\|{\mathbf Z}_1 {\mathbf b}_1\|_2 +
\cdots
+ \|{\mathbf Z}_q{\mathbf b}_q\|_2.
\]
Since ${\mathbf Z}_A {\mathbf b}= {\mathbf Z}_1 {\mathbf b}_1 + \cdots
+ {\mathbf Z}_q {\mathbf b}_q$, the above two
inequalities imply that
\[
c_3 (\|{\mathbf Z}_1 {\mathbf b}_1\|_2 + \cdots+ \|{\mathbf
Z}_q{\mathbf b}_q\|_2)
\le\|{\mathbf Z}_A {\mathbf b}\|_2 \le\|{\mathbf Z}_1 {\mathbf b}_1\|
_2 + \cdots+ \|{\mathbf Z}_q{\mathbf b}
_q\|_2.
\]
Therefore,
%
%
\begin{eqnarray}
\label{Sa}
&&
c_3^2(\|{\mathbf Z}_1 {\mathbf b}_1\|_2^2 + \cdots+ \|{\mathbf
Z}_q{\mathbf b}_q\|_2^2)\nonumber\\[-8pt]\\[-8pt]
&&\qquad
\le\|{\mathbf Z}_A {\mathbf b}\|_2^2 \le2 (\|{\mathbf Z}_1 {\mathbf
b}_1\|_2^2 + \cdots+ \|{\mathbf Z}
_q{\mathbf b}
_q\|_2^2).\nonumber
\end{eqnarray}
Let ${\mathbf C}_j = n^{-1}{\mathbf Z}_j'{\mathbf Z}_j$. By Lemma 6.2
of Zhou, Shen and
Wolf (\citeyear{ZShW98}),
%
%
\begin{equation}
\label{Sb}
c_4h \le\rho_{\min}({\mathbf C}_j) \le\rho_{\max}({\mathbf C}_j)
\le c_5 h,\qquad j \in A.
\end{equation}
Since ${\mathbf C}_A = n^{-1} {\mathbf Z}_A' {\mathbf Z}_A$, it follows
from (\ref{Sa}) that
\[
c_3^2 ({{\mathbf b}_1'{\mathbf C}_1{\mathbf b}_1}
+ \cdots+
{{\mathbf b}_q'{\mathbf C}_q{\mathbf b}_q} )
\le{{\mathbf b}'{\mathbf C}_A{\mathbf b}} \le
2 ({{\mathbf b}_1'{\mathbf C}_1{\mathbf b}_1} + \cdots+
{{\mathbf b}_q'{\mathbf C}_q{\mathbf b}_q} ).
\]
Therefore, by (\ref{Sb}),
\begin{eqnarray*}
\frac{{\mathbf b}_1'{\mathbf C}_1{\mathbf b}_1}{\|{\mathbf b}\|_2^2}
+ \cdots+
\frac{{\mathbf b}_q'{\mathbf C}_q{\mathbf b}_q}{\|{\mathbf b}\|_2^2}
&=&
\frac{{\mathbf b}_1'{\mathbf C}_1{\mathbf b}_1}{\|{\mathbf b}_1\|
_2^2}\frac{\|{\mathbf b}_1\|_2^2}{\|{\mathbf b}
\|_2^2}
+ \cdots+
\frac{{\mathbf b}_q'{\mathbf C}_q{\mathbf b}_q}{\|{\mathbf b}_q\|
_2^2}\frac{\|{\mathbf b}_q\|_2^2}{\|{\mathbf b}
\|_2^2}
\\
& \ge& \rho_{\min}({\mathbf C}_1)\frac{\|{\mathbf b}_1\|_2^2}{\|
{\mathbf b}\|_2^2}
+ \cdots+ \rho_{\min}({\mathbf C}_q)\frac{\|{\mathbf b}_q\|_2^2}{\|
{\mathbf b}\|_2^2} \\
& \ge& c_4 h.
\end{eqnarray*}
Similarly,
\[
\frac{{\mathbf b}_1'{\mathbf C}_1{\mathbf b}_1}{\|{\mathbf b}\|_2^2}
+ \cdots+
\frac{{\mathbf b}_q'{\mathbf C}_q{\mathbf b}_q}{\|{\mathbf b}\|_2^2}
\le c_5 h.
\]
Thus, we have
\[
c_3^2 c_4 h \le\frac{{\mathbf b}'{\mathbf C}_A {\mathbf b}}{{\mathbf
b}'{\mathbf b}} \le2 c_5 h.
\]
The lemma follows.
\end{pf}
\begin{pf*}{Proof of Theorem \ref{ThmA}}
The proof of parts (i) and (ii) essentially follows the
proof of Theorem 2.1 of Wei and Huang (\citeyear{WH08}). The only change that must
be made here is that we need to consider the approximation error of the
regression functions by splines. Specifically, let $\bolds\xi
_n=\bolds\varepsilon_n +
\bolds\delta_n$, where
$\bolds\delta_n=(\delta_{n1}, \ldots, \delta_{nn})'$ with $\delta
_{ni}=\sum
_{j=1}^{q_n}(f_{0j}(X_{ij})-f_{nj}(X_{ij}))$.
Since $\|f_{0j}-f_{nj}\|_{2}=O(m_n^{-d})=O(n^{-d/(2d+1)})$ for
$m_n=n^{1/(2d+1)}$,
we have
\[
\|\bolds\delta_n\|_2
\le C_1 \sqrt{n q m_n^{-2d}}= C_1 q n^{1/(4d+2)}
\]
for some constant $C_1 > 0$.
For any integer $t$, let
\[
\chi_t= \max_{|A|=t} \max_{\|U_{A_k}\|_2=1, 1\le k \le t}
\frac{|\bolds\xi_n'V_A({\mathbf s})|}{\|V_A({\mathbf s})\|_2}
\quad\mbox{and}\quad
\chi_t^*= \max_{|A|=t} \max_{\|U_{A_k}\|_2=1, 1\le k \le t}
\frac{|\bolds\varepsilon_n'V_A({\mathbf s})|}{\|V_A({\mathbf s})\|_2},
\]
where $V_A(S_A)= \bolds\xi_n^{\prime}({\mathbf Z}_{A}({\mathbf
Z}_{A}^{\prime}{\mathbf Z}_{A})^{-1}
\bar{S}_{A}-(I-P_{A})X\bolds\beta$
for $N(A)=q_{1}=m\geq0$, ${S}_{A}=({S}_{A_{1}}^{\prime},\ldots
,{S}_{A_{m}}^{\prime})^{\prime}$,
${S}_{A_{k}}=\lambda\sqrt{d_{A_{k}}}U_{A_{k}}$ and $\|U_{A_{k}}\|_{2}=1$.

For a sufficiently large constant $C_2 > 0$, define
\[
\Omega_{t_{0}}=\bigl\{({\mathbf Z}, \bolds\varepsilon_n)\dvtx x_{t} \leq
\sigma C_2
\sqrt{(t \vee1)m_n\log(pm_n)},\forall t\geq t_{0}\bigr\}
\]
and
\[
\Omega_{t_{0}}^*=\bigl\{({\mathbf Z}, \bolds\varepsilon_n)\dvtx
x_{t}^{*}\leq\sigma C_2
\sqrt{(t \vee1)m_n\log(p m_n)},\forall t\geq t_{0}\bigr\},
\]
where $t_0 \ge0$.

As in the proof of Theorem 2.1 of Wei and Huang (\citeyear{WH08}),
\[
({\mathbf Z}, \bolds\varepsilon_n) \in\Omega_q \quad\Rightarrow\quad
|{\widetilde A}_1| \le M_1 q
\]
for a constant $M_1 > 1$.
By the triangle and Cauchy--Schwarz inequalities,
%
%
\begin{equation}
\label{Oa}
\frac{|\bolds\xi_n'V_A({\mathbf s})|}{\|V_A({\mathbf s})\|_2}
=\frac{|\bolds\varepsilon_n'V_A({\mathbf s}
)+\bolds\delta_n'V_A({\mathbf s})|}{\|V_A({\mathbf s})\|_2} \le
\frac{|\bolds\varepsilon_n'V_A({\mathbf s}
)|}{\|
V_A\|_2}+\|\bolds\delta_n\|.
\end{equation}
In the proof of Theorem 2.1 of Wei and Huang (\citeyear{WH08}), it is shown that
%
%
\begin{equation}
\label{Ob}
\mathrm{P}(\Omega_0^*) \ge2- \frac{2}{p^{1+c_0}}-\exp\Bigl(\frac
{2p}{p^{1+c_0}} \Bigr) \rightarrow1.
\end{equation}
Since
\[
\frac{|\bolds\delta_n'V_A({\mathbf s})|}{\|V_A({\mathbf s})\|_2}
\le\|\bolds\delta_n\|_2
\le C_1 {q} n^{{1}/({2(2d+1)})}
\]
and $m_n=O(n^{1/(2d+1)})$, we have for all $t \ge0$ and $n$
sufficiently large,
%
%
\begin{equation}
\label{Oc}
\|\bolds\delta_n\|_2 \le C_1 {q} n^{{1}/({2(2d+1)})}\le\sigma
C_2 \sqrt
{(t\vee1)m_n\log(p)} .
\end{equation}
It follows from (\ref{Oa}), (\ref{Ob}) and (\ref{Oc}) that
\(
\mathrm{P}(\Omega_0) \rightarrow1.
\)
This completes the proof of part (i) of Theorem \ref{ThmA}.

Before proving part (ii), we first prove part (iii) of Theorem
\ref{ThmA}. By the definition of
$\widetilde{\bolds\beta}_n \equiv(\widetilde{\bolds\beta}_{n1}',
\ldots, \widetilde{\bolds\beta}_{np}')'$,
%
%
\begin{equation}
\label{DefA} \|{\mathbf Y}-{\mathbf Z}\widetilde{\bolds\beta}_n\|
_2^2+ \lambda_{n1} \sum_{j=1}^{p}
\|\widetilde{\bolds\beta}_{nj}\|_2 \le\|{\mathbf Y}-{\mathbf
Z}\bolds\beta_n\|_2^2 + \lambda_{n1}
\sum_{j=1}^{p} \|\bolds\beta_{nj}\|_2.
\end{equation}
Let $A_2=\{j\dvtx\|\bolds\beta_{nj}\|_2 \neq0 \mbox{ or } \|
\widetilde{\bolds\beta}
_{nj}\|_2
\neq0\}$ and $d_{n2}=|A_2|$. By part (i), $d_{n2}=O_p(q)$.
By (\ref{DefA}) and the definition of $A_2$,
%
%
\begin{eqnarray}
\label{DefB}
&&
\|{\mathbf Y}-{\mathbf Z}_{A_2}\widetilde{\bolds\beta
}_{nA_2}\|_2^2+ \lambda_{n1}
\sum_{j
\in A_2}\|\widetilde{\bolds\beta}_{nj}\|_2 \nonumber\\[-8pt]\\[-8pt]
&&\qquad\le\|{\mathbf
Y}-{\mathbf Z}_{A_2}\bolds\beta_{nA_2}\|_2^2 +
\lambda_{n1} \sum_{j\in A_2}\|\bolds\beta_{nj}\|_2.\nonumber
\end{eqnarray}
Let $\bolds\eta_n={\mathbf Y}-{\mathbf Z}\bolds\beta_n$. Write
\[
{\mathbf Y}-{\mathbf Z}_{A_2}\widetilde{\bolds\beta}_{nA_2} =
{\mathbf Y}-{\mathbf Z}\bolds\beta_n -
{\mathbf Z}_{A_2}(\widetilde{\bolds\beta}_{nA_2}-\bolds\beta_{nA_2})
=\bolds\eta_n-{\mathbf Z}_{A_2}(\widetilde{\bolds\beta
}_{nA_2}-\bolds\beta_{nA_2}).
\]
We have
\[
\|{\mathbf Y}-{\mathbf Z}_{A_2}\widetilde{\bolds\beta}_{nA_2}\|
_2^2=\|{\mathbf Z}_{A_2}(\widetilde{\bolds\beta}_{nA_2}-
\bolds\beta_{nA_2})\|_2^2- 2 \bolds\eta_n'{\mathbf
Z}_{A_2}(\widetilde{\bolds\beta}_{nA_2}-\bolds\beta_{nA_2})
+ \bolds\eta_n'\bolds\eta_n.
\]
We can rewrite (\ref{DefB}) as
%
%
\begin{eqnarray}
\label{DefC}
&&\|{\mathbf Z}_{A_2}(\widetilde{\bolds\beta}_{nA_2}-
\bolds\beta_{nA_2})\|_2^2 - 2 \bolds\eta_n'{\mathbf
Z}_{A_2}(\widetilde{\bolds\beta}_{nA_2}-\bolds\beta_{nA_2})
\nonumber\\[-8pt]\\[-8pt]
&&\qquad
\le\lambda_{n1} \sum_{j\in A_1} \|\bolds\beta_{nj}\|_2 -\lambda_{n1}
\sum_{j\in A_1} \|\widetilde{\bolds\beta}_{nj}\|_2.\nonumber
\end{eqnarray}
Now
%
%
\begin{eqnarray}
\label{DefD}
\biggl| \sum_{j\in A_1} \|\bolds\beta_{nj}\|_2-\sum_{j\in
A_1}\|\widetilde{\bolds\beta}_{nj}\|_2 \biggr| &\le& \sqrt{|A_1|} \cdot
\|\widetilde{\bolds\beta}_{nA_1}-\bolds\beta_{nA_1}\|_2 \nonumber\\[-8pt]\\[-8pt]
&\le& \sqrt
{|A_1|} \cdot
\|\widetilde{\bolds\beta}_{nA_2}-\bolds\beta_{nA_2}\|_2.\nonumber
\end{eqnarray}
Let $\bolds\nu_n ={\mathbf Z}_{A_2}(\widetilde{\bolds\beta
}_{nA_2}-\bolds\beta_{nA_2})$. Combining
(\ref{DefB}), (\ref{DefC}) and (\ref{DefD}) to get
%
%
\begin{equation}
\label{DefE} \|\bolds\nu_n\|_2^2 - 2 \bolds\eta_n'\bolds\nu_n
\le\lambda_{n1}
\sqrt{|A_1|} \cdot\|\widetilde{\bolds\beta}_{nA_2}-\bolds\beta
_{nA_2}\|_2.
\end{equation}

Let $\bolds\eta_n^*$ be the projection of $\bolds\eta_n$ to the
span of
${\mathbf Z}_{A_2}$, that is, $\bolds\eta_n^* =
{\mathbf Z}_{A_2}({\mathbf Z}_{A_2}'\times\break{\mathbf Z}_{A_2})^{-1}{\mathbf
Z}_{A_2}'\bolds\eta_n$. By the
Cauchy--Schwarz inequality,
%
%
\begin{equation}
\label{DefF} 2|\bolds\eta_n'\bolds\nu_n| \le2 \|\bolds\eta
_n^*\|_2\cdot\|\bolds\nu
_n\|_2
\le2 \|\bolds\eta_n^*\|_2^2 + \tfrac{1}{2}\|\bolds\nu_n\|_2^2.
\end{equation}
From (\ref{DefE}) and (\ref{DefF}), we have
\[
\|\bolds\nu_n\|_2^2 \le4 \|\bolds\eta_n^*\|_2^2 + 2 \lambda
_{n1} \sqrt{|A_1|}
\cdot\|\widetilde{\bolds\beta}_{nA_2}-\bolds\beta_{nA_2}\|_2.
\]
Let $c_{n*}$ be the smallest eigenvalue of
${\mathbf Z}_{A_2}'{\mathbf Z}_{A_2}/n$. By Lemma \ref{LemPB} and part (i),
$c_{n*}\asymp_p m_n^{-1}$.
Since $\|\bolds\nu_n\|_2^2 \ge n c_{n*}
\|\widetilde{\bolds\beta}_{nA_2}-\bolds\beta_{nA_2}\|_2^2$ and
$2ab \le a^2 + b^2$,
\[
nc_{n*}\|\widetilde{\bolds\beta}_{nA_2}-\bolds\beta_{nA_2}\|_2^2
\le4 \|\bolds\eta_n^*\|_2^2 +
\frac{(2\lambda_{n1}
\sqrt{|A_1|})^2}{2nc_{n*}}+\frac{1}{2}nc_{n*}\|\widetilde{\bolds
\beta}_{nA_2}-\bolds\beta
_{nA_2}\|_2^2.
\]
It follows that
%
%
\begin{equation}
\label{Da} \|\widetilde{\bolds\beta}_{nA_2}-\bolds\beta_{nA_2}\|
_2^2 \le\frac{8
\|\bolds\eta_n^*\|_2^2}{nc_{n*}} + \frac{4 \lambda_{n1}^2
|A_1|}{n^2c_{n*}^2}.
\end{equation}
Let $f_0({\mathbf X}_i) = \sum_{j=1}^{p}f_{0j}(X_{ij})$ and
$f_{0A}({\mathbf X}_i)=\sum_{j \in A} f_{0j}(X_{ij})$. Write
\begin{eqnarray*}
\eta_i &=& Y_i-\mu- f_0({\mathbf X}_i) + (\mu-\overline{Y}) +
f_0({\mathbf X}_i)-\sum_{j\in A_2} Z_{ij}'\bolds\beta
_{nj}\\
&=&\varepsilon_i +(\mu-\overline{Y})
+ f_{A_2}({\mathbf X}_i)-f_{nA_2}({\mathbf X}_i).
\end{eqnarray*}
Since $|\mu-\overline{Y}|^2=O_p(n^{-1})$ and
$\|f_{0j}-f_{nj}\|_{\infty}=O(m_n^{-d})$, we have
%
%
\begin{equation}
\label{Db}
\|\bolds\eta_n^*\|_2^2 \le2\|\bolds\varepsilon_n^*\|_2^2 +
O_p(1) + O(n d_{n2}
m_n^{-2d}),
\end{equation}
where $\bolds\varepsilon_n^*$ is the projection of $\bolds
\varepsilon_n=(\varepsilon_1,
\ldots, \varepsilon_n)'$ to the span of ${\mathbf Z}_{A_2}$. We have
\[
\|\bolds\varepsilon_n^*\|_2^2 =
\|({\mathbf Z}_{A_2}'{\mathbf Z}_{A_2})^{-1/2}{\mathbf Z}_{A_2}'\bolds
\varepsilon_n\|_2^2 \le
\frac{1}{nc_{n*}} \|{\mathbf Z}_{A_2}'\bolds\varepsilon_n\|_2^2.
\]
Now
\[
{\max_{A\dvtx|A| \le d_{n2}}}\|{\mathbf Z}_{A}'\bolds\varepsilon_n\|
_2^2 = \max_{A\dvtx
|A| \le
d_{n2}} \sum_{j \in A}\|{\mathbf Z}_j'\bolds\varepsilon_n\|_2^2 \le
{d_{n2} m_n \max
_{1\le
j\le p, 1 \le k \le m_n}}|\mathcal{Z}_{jk}'\bolds\varepsilon|^2,
\]
where $\mathcal{Z}_{jk} = (\psi_{k}(X_{1j}), \ldots,
\psi_{k}(X_{nj}))'$.
By Lemma \ref{LemC},
\begin{eqnarray*}
{\max_{1\le j\le p, 1\le k \le m_n}}|\mathcal{Z}_{jk}'\bolds
\varepsilon_n|^2 &=& {n m_n^{-1}
\max_{1\le j\le p, 1 \le k \le
m_n}}|(m_n/n)^{1/2}\mathcal{Z}_{jk}'\bolds\varepsilon_n|^2\\
&=& O_p(1) n
m_n^{-1}\log(p
m_n).
\end{eqnarray*}
It follows that,
%
%
\begin{equation}
\label{Dc} \|\bolds\varepsilon_n^*\|_2^2 = O_p(1)\frac{d_{n2} \log(p
m_n)}{c_{n*}}.
\end{equation}
Combining (\ref{Da}), (\ref{Db}) and (\ref{Dc}), we get
\begin{eqnarray*}
\|\widetilde{\bolds\beta}_{A_2}-\bolds\beta_{A_2}\|_2^2 &\le& O_p
\biggl(\frac{d_{n2}
\log(p m_n)}{nc_{n*}^2} \biggr) +
O_p \biggl(\frac{1}{nc_{n*}} \biggr)\\
&&{}+O \biggl(\frac{d_{n2}
m_n^{-2d}}{c_{n*}} \biggr) + \frac{4 \lambda_{n1}^2 |A_1|}{n^2c_{n*}^2}.
\end{eqnarray*}
Since $d_{n2}=O_p(q)$, $c_{n*}\asymp_p m_n^{-1}$ and $c_n^{*}\asymp_p
m_n^{-1}$, we have
\[
\|\widetilde{\bolds\beta}_{A_2}-\bolds\beta_{A_2}\|_2^2 \le
O_p \biggl(\frac{m_n^2 \log(p m_n)}{n} \biggr) +
O_p \biggl(\frac{m_n}{n} \biggr)+O \biggl(\frac{1}{
m_n^{2d-1}} \biggr) + O \biggl(\frac{4 m_n^2\lambda_{n1}^2
}{n^2} \biggr).
\]
This completes the proof of part (iii).

We now prove part (ii). Since $\|f_j\|_2 \ge c_f>0, 1\le j \le q$, $\|
f_j-f_{nj}\|_2 = O(m_n^{-d})$
and $\|f_{nj}\|_2 \ge\|f_j\|_2 - \|f_j-f_{nj}\|_2$, we have
$\|f_{nj}\|_2 \ge0.5 c_f $ for $n$ sufficiently large.
By a result of de Boor (\citeyear{deBoor01}), see also (12) of Stone (\citeyear{Stone86}), there are
positive constants $c_6$ and $c_7$ such that
\[
c_6 m_n^{-1}\|\bolds\beta_n\|_2^2 \le\|f_{nj}\|_2^2\le c_7 m_n^{-1}
\|
\bolds\beta
_{nj}\|_2^2.
\]
It follows that
\(
\|\bolds\beta_{nj}\|_2^2 \ge c_7^{-1} m_n \|f_{nj}\|_2^2 \ge
0.25 c_7^{-1} c_f^2 m_n.
\)
Therefore, if $\|\bolds\beta_{nj}\|_2 \neq0$ but $\|\widetilde
{\bolds\beta}_{nj}\|_2=0$,
then
%
%
\begin{equation}
\label{ContrA}
\|\widetilde{\bolds\beta}_{nj}-\bolds\beta_{nj}\|_2^2 \ge0.25
c_7^{-1} c_f^2 m_n.
\end{equation}
However, since $(m_n\log(pm_n))/n \rightarrow0$ and
$(\lambda_{n1}^2 m_n)/n^2 \rightarrow$, (\ref{ContrA}) contradicts
part (iii).
\end{pf*}
\begin{pf*}{Proof of Theorem \ref{ThmB}}
By the definition of $\widetilde{f}_j, 1 \le j \le p$, parts (i) and
(ii) follow
from parts (i) and (ii) of Theorem \ref{ThmA} directly.

Now consider part (iii). By the properties of spline [de Boor (\citeyear{deBoor01})],
\[
c_6 m_n^{-1} \|\widetilde{\bolds\beta}_{nj}-\bolds\beta_{nj}\|_2^2
\le\|\widetilde{f}_{nj}-f_{nj}\|_2^2 \le
c_7 m_n^{-1} \|\widetilde{\bolds\beta}_{nj}-\bolds\beta_{nj}\|_2^2.
\]
Thus,
%
%
\begin{eqnarray}
\label{SSa}
\|\widetilde{f}_{nj}-f_{nj}\|_2^2 &=& O_p \biggl(\frac{ m_n \log(p m_n)}{n} \biggr)+
O_p \biggl(\frac{1}{n} \biggr)\nonumber\\[-8pt]\\[-8pt]
&&{} +
O \biggl(\frac{1}{
m_n^{2d}} \biggr) + O \biggl(\frac{4 m_n\lambda_{n1}^2}{n^2}
\biggr).\nonumber
\end{eqnarray}
By (A3),
%
%
\begin{equation}
\label{SSb}
\|f_j-f_{nj}\|_2^2 = O(m_n^{-2d}).
\end{equation}
Part (iii) follows from (\ref{SSa}) and (\ref{SSb}).
\end{pf*}

In the proofs below, for any matrix ${\mathbf H}$, denote
its $2$-norm by $\|{\mathbf H}\|$, which is equal to its largest eigenvalue.
This norm satisfies the inequality $\|{\mathbf H}{\mathbf x}\| \le\|
{\mathbf H}\| \|
{\mathbf x}\|$
for a
column vector ${\mathbf x}$ whose dimension is the same as the number of
the columns of ${\mathbf H}$.

Denote $\bolds\beta_{nA_1} = (\bolds\beta_{nj}', j \in A_1)'$,
$\widehat{\bolds\beta}_{nA_1}=(\widehat{\bolds\beta}_{nj}', j\in
A_1)'$ and ${\mathbf Z}_{A_1}=({\mathbf Z}_j, j
\in A_1)$. Define ${\mathbf C}_{A_1}=n^{-1}{\mathbf Z}_{A_1}'{\mathbf Z}_{A_1}$.
Let $\rho_{n1}$ and $\rho_{n2}$ be the smallest and largest eigenvalues
of ${\mathbf C}_{A_1}$, respectively.
\begin{pf*}{Proof of Theorem \ref{ThmC}}
By the KKT, a necessary and sufficient condition for $\widehat{\bolds
\beta}_n$
is
%
%
\begin{equation}
\label{KKTb} \cases{
2{\mathbf Z}_j' ({\mathbf Y}-{\mathbf Z}\widehat{\bolds\beta}_n )
= \lambda_{n2} w_{nj} \dfrac{\widehat{\bolds\beta}_{nj}}{\|\widehat
{\bolds\beta}_{nj}\|}, &\quad
$\|\widehat{\bolds\beta}_j\|_2 \neq0, j \ge1$, \vspace*{2pt}\cr
2\|{\mathbf Z}_j' ({\mathbf Y}-{\mathbf Z}\widehat{\bolds\beta}_n )\|
_2 \le\lambda_{n2} w_{nj}, &\quad
$\|\widehat{\bolds\beta}_{nj}\|=0, j \ge1$.}
\end{equation}

Let $\bolds\nu_n=(w_{nj} \widehat{\bolds\beta}_j/(2\|\widehat
{\bolds\beta}_{nj}\|), j \in A_1)'$. Define
%
%
\begin{equation}
\label{OracleA} \widehat{\bolds\beta}_{nA_1} =
({\mathbf Z}_{A_1}'{\mathbf Z}_{A_1})^{-1}({\mathbf Z}_{A_1}'{\mathbf
Y}- \lambda_{n2} \bolds\nu_n).
\end{equation}

If $\widehat{\bolds\beta}_{nA_1} =_0 \bolds\beta_{nA_1}$, then
the equation in
(\ref{KKTb}) holds for $\widehat{\bolds\beta}_n \equiv(\widehat
{\bolds\beta}_{nA_1}',
{\mathbf0}')'$. Thus, since ${\mathbf Z}\widehat{\bolds\beta
}_n={\mathbf Z}_{A_1}\widehat{\bolds\beta}_{nA_1}$ for
this $\widehat{\bolds\beta}_n$ and $\{{\mathbf Z}_{j},j\in A_1\}$
are linearly
independent,
\[
\widehat{\bolds\beta}_n =_0 \bolds\beta_{n} \qquad\mbox{if }
\cases{
\widehat{\bolds\beta}_{nA_1}=_0 \bolds\beta_{nA_1}, \vspace*{2pt}\cr
\|{\mathbf Z}_{j}' ({\mathbf Y}-{\mathbf Z}_{A_1}\widehat{\bolds\beta
}_{nA_1} )\|_2 \le
\lambda_{n2} w_{nj}/2, &\quad $\forall j\notin A_1$.}
\]
This is true if
\[
\widehat{\bolds\beta}_n =_0 \bolds\beta_{n} \qquad\mbox{if }
\cases{
\|\bolds\beta_{nj}\|_2 -\|\widehat{\bolds\beta}_{nj}\|_2 < \|
\bolds\beta_{nj}\|_2, &\quad
$\forall j \in A_1$, \vspace*{2pt}\cr
\|{\mathbf Z}_{j}'({\mathbf Y}-{\mathbf Z}_{A_1}\widehat{\bolds\beta
}_{nA_1})\|_2 \le
\lambda_{n2} w_{nj}/2, &\quad $\forall j\notin A_1$.}
\]
Therefore,
\begin{eqnarray*}
\mathrm{P}(\widehat{\bolds\beta}_n \neq_0 \bolds\beta_n) &\le&
\mathrm{P}(\|\widehat{\bolds\beta}_{nj}-\bolds\beta_{nj}\|_2 \ge
\|\bolds\beta_{nj}\|_2,
\exists
j \in A_1 ) \\
&&{} +
\mathrm{P}\bigl(\|{\mathbf Z}_{j}'({\mathbf Y}-{\mathbf Z}_{A_1}\widehat
{\bolds\beta}_{nA_1})\|_2
> \lambda_{n2} w_{nj}/2, \exists j \notin A_1 \bigr).
\end{eqnarray*}

Let $f_{0j}({\mathbf X}_j) = (f_{0j}(X_{1j}), \ldots,
f_{0j}(X_{nj}))'$ and $\bolds\delta_n = \sum_{j \in A_1}
f_{0j}({\mathbf X}_j)-\break {\mathbf Z}_{A_1}\bolds\beta_{nA_1}$. By Lemma
\ref{LemA}, we
have
%
%
\begin{equation}
\label{DeltaA} n^{-1}\|\bolds\delta_n\|^2 =O_p( q m_n^{-2d}).
\end{equation}
Let ${\mathbf H}_n={\mathbf I}_n-{\mathbf Z}_{A_1}({\mathbf
Z}_{A_1}'{\mathbf Z}_{A_1})^{-1}{\mathbf Z}_{A_1}'$. By
(\ref{OracleA}),
%
%
\begin{equation}
\label{EEa} \widehat{\bolds\beta}_{nA_1}-\bolds\beta_{nA_1} = n^{-1}
{\mathbf C}_{A_1}^{-1} \bigl({\mathbf Z}_{A_1}'(\bolds\varepsilon
_n+\bolds\delta_n)-\lambda_{n2}
\bolds\nu_n \bigr)
\end{equation}
and
%
%
\begin{equation}
\label{EEb} {\mathbf Y}-{\mathbf Z}_{A_1}\widehat{\bolds\beta
}_{nA_1} = {\mathbf H}_n \bolds\varepsilon_n + {\mathbf H}_n
\bolds\delta_n+ \lambda_{n2} {\mathbf Z}_{A_1}{\mathbf
C}_{A_1}^{-1}\bolds\nu_n/n.
\end{equation}
Based on these two equations, Lemma \ref{LemE} below shows that
\[
\mathrm{P}( \|\widehat{\bolds\beta}_{nj} - \bolds\beta_{nj}\|_2
\ge\|\bolds\beta_{nj}\|_2,
\exists j \in A_1 ) \rightarrow0,
\]
and Lemma \ref{LemF} below shows that
\[
\mathrm{P}\bigl(
\|{\mathbf Z}_{j}'({\mathbf Y}-{\mathbf Z}_{A_1}\widehat{\bolds\beta
}_{nA_1})\|_2 >
\lambda_{n2} w_{nj}/2, \exists j\notin A_1 \bigr) \rightarrow0.
\]
These two equations lead to part (i) of the theorem.

We now prove part (ii) of Theorem \ref{ThmC}. As in (\ref{Db}), for
$\bolds\eta_n={\mathbf Y}-{\mathbf Z}\bolds\beta_n$ and
\[
\bolds\eta_{n1}^* = {\mathbf Z}_{A_1}({\mathbf Z}_{A_1}'{\mathbf
Z}_{A_1})^{-1}{\mathbf Z}_{A_1}'\bolds\eta_{n},
\]
we have
%
%
\begin{equation}
\label{Ca}
\|\bolds\eta_{n1}^*\|_2^2 \le2\|\bolds\varepsilon_{n1}^*\|_2^2 +
O_p(1) + O(q n
m_n^{-2d}),
\end{equation}
where $\bolds\varepsilon_{n1}^*$ is the projection of $\bolds
\varepsilon_n=(\varepsilon_1,
\ldots, \varepsilon_n)'$ to the span of ${\mathbf Z}_{A_1}$. We have
%
%
\begin{equation}
\label{Cb}
\|\bolds\varepsilon_{n1}^*\|_2^2 =
\|({\mathbf Z}_{A_1}'{\mathbf Z}_{A_1})^{-1/2}{\mathbf Z}_{A_1}'\bolds
\varepsilon_n\|_2^2 \le
\frac{1}{n\rho_{n1}} \|{\mathbf Z}_{A_1}'\bolds\varepsilon_n\|_2^2
= O_p(1) \frac{|A_1|}{\rho_{n1}}.
\end{equation}
Now similarly to the proof of (\ref{Da}), we can show that
%
%
\begin{equation}
\label{Cc} \|\widehat{\bolds\beta}_{nA_1}-\bolds\beta_{nA_1}\|
_2^2 \le\frac{8
\|\bolds\eta_{n1}^*\|_2^2}{n\rho_{n1}} + \frac{4 \lambda_{n2}^2
|A_1|}{n^2\rho_{n1}^2}.
\end{equation}
Combining (\ref{Ca}), (\ref{Cb}) and (\ref{Cc}), we get
\[
\|\widehat{\bolds\beta}_{nA_1}-\bolds\beta_{nA_1}\|_2^2=
O_p \biggl(\frac{8 }{n\rho_{n1}^2} \biggr) +
O_p \biggl(\frac{1}{n \rho_{n1}} \biggr)+O \biggl(\frac{1}{
m_n^{2d-1}} \biggr) + O \biggl(\frac{4 \lambda_{n2}^2
}{n^2\rho_{n1}^2} \biggr).
\]
Since $\rho_{n1}\asymp_p m_n^{-1}$, the result follows.
\end{pf*}

The following lemmas are needed in the proof of Theorem \ref{ThmC}.
\begin{lemma}
\label{LemD} For $\bolds\nu_n=(w_{nj} \widetilde{\bolds\beta
}_j/(2\|\widetilde{\bolds\beta}_{nj}\|), j
\in A_1)'$, under condition \textup{(B1)},
\[
\|\bolds\nu_n\|^2 =O_p(h_n^2)=
O_p \bigl((b_{n1}^2c_b)^{-2}r_n^{-1}+qb_{n1}^{-1} \bigr).
\]
\end{lemma}
\begin{pf}
Write
\[
\|\bolds\nu_n\|^2 = \sum_{j\in A_1} w_{j}^2 = \sum_{j\in A_1}
\|\widetilde{\bolds\beta}_{nj}\|^{-2} = \sum_{j\in A_1}
\frac{\|\bolds\beta_{nj}\|^2-\|\widetilde{\bolds\beta}_{nj}\|^2}
{\|\bolds\beta_{nj}\|^2\cdot\|\widetilde{\bolds\beta}_{nj}\|^2}
+ \sum_{j\in
A_1}\|\bolds\beta_{nj}\|^{-1}.
\]
Under (B2),
\[
\sum_{j\in A_1}
\frac{ |\|\bolds\beta_{nj}\|^2-\|\widetilde{\bolds\beta}_{nj}\|
^2 |}
{\|\bolds\beta_{nj}\|^2\cdot\|\widetilde{\bolds\beta}_{nj}\|^2}
\le M c_b^{-2}
b_{n1}^{-4} \|\widetilde{\bolds\beta}_{n}-\bolds\beta_{n}\|
\]
and
\(
\sum_{j\in A_1}\|\bolds\beta_{nj}\|^{-2} \le q b_{n1}^{-2}.
\)
The claim follows.
\end{pf}

Let $\rho_{n3}$ be the maximum of the
largest eigenvalues of $n^{-1}{\mathbf Z}_j'{\mathbf Z}_j, j \in A_0$, that
is, $\rho_{n3}=\max_{j \in A_0} \|n^{-1}{\mathbf Z}_j'{\mathbf Z}_j\|_2$.
By Lemma \ref{LemPB},
%
%
\begin{equation}
\label{Qnote}
b_{n1} \asymp O(m_n^{1/2}),\qquad
\rho_{n1}\asymp_p m_n^{-1}, \qquad\rho_{n2}\asymp_p m_n^{-1}
\quad\mbox{and}\quad
\rho_{n3}\asymp_p m_n^{-1}.\hspace*{-27pt}
\end{equation}
\begin{lemma}
\label{LemE} Under conditions \textup{(B1)}, \textup{(B2)}, \textup{(A3)} and
\textup{(A4)},
%
%
\begin{equation}
\label{KKTaa} \mathrm{P}( \|\widehat{\bolds\beta}_{nj} - \bolds
\beta_{nj}\|_2 \ge
\|\bolds\beta_{nj}\|_2, \exists j \in A_1 ) \rightarrow0.
\end{equation}
\end{lemma}
\begin{pf}
Let ${\mathbf T}_{nj}$ be an $m_n \times
q m_n$ matrix with the form
\[
{\mathbf T}_{nj} = ({\mathbf0}_{m_n}, \ldots, {\mathbf0}_{m_n},
{\mathbf I}_{m_n},
{\mathbf0}_{m_n}, \ldots, {\mathbf0}_{m_n}),
\]
where ${\mathbf O}_{m_n}$ is an $m_n \times m_n$ matrix of zeros and
${\mathbf I}_{m_n}$ is an $m_n \times m_n$ identity matrix, and
${\mathbf I}_{m_n}$ is at the $j$th block. By (\ref{EEa}),
\(
\widehat{\bolds\beta}_{nj} - \bolds\beta_{nj} = n^{-1}{\mathbf T}_{nj}
{\mathbf C}_{A_1}^{-1}({\mathbf Z}_{A_1}'\bolds\varepsilon_n +
{\mathbf Z}_{A_1}'\bolds\delta_n-\lambda_{n2}
\bolds\nu_n).
\)
By the triangle inequality,
%
%
\begin{eqnarray}
\label{TrA}
\|\widehat{\bolds\beta}_{nj} - \bolds\beta_{nj}\|_2 &\le& n^{-1} \|
{\mathbf T}_{nj}
{\mathbf C}_{A_1}^{-1}{\mathbf Z}_{A_1}'\bolds\varepsilon_n\|_2 +
n^{-1}\|{\mathbf T}_{nj}
{\mathbf C}_{A_1}^{-1}{\mathbf
Z}_{A_1}'\delta_n\|_2\nonumber\\[-8pt]\\[-8pt]
&&{} + n^{-1}\lambda
_{n2} \|{\mathbf T}_{nj}
{\mathbf C}_{A_1}^{-1}\bolds\nu_n\|_2.\nonumber
\end{eqnarray}
Let $C$ be a generic constant independent of $n$.
The first term on the right-hand side
%
%
\begin{eqnarray}
\label{Ta}
{\max_{j \in A_1} n^{-1} }\|{\mathbf T}_{nj} {\mathbf
C}_{A_1}^{-1}{\mathbf Z}_{A_1}'\bolds\varepsilon
_n\|_2
&\le&
n^{-1} \rho_{n1}^{-1} \|{\mathbf Z}_{A_1}'\bolds\varepsilon_n\|_2
\nonumber\\
&=& n^{-1/2} \rho_{n1}^{-1} \|n^{-1/2}{\mathbf Z}_{A_1}'\bolds
\varepsilon_n\|_2
\\
&= & O_p(1)
n^{-1/2}\rho_{n1}^{-1}m_n^{-1/2}(q m_n)^{1/2}.\nonumber
\end{eqnarray}
By (\ref{DeltaA}), the second term
%
%
\begin{eqnarray}
\label{Tb} \max_{j \in A_1} n^{-1}\|{\mathbf T}_{nj}
{\mathbf C}_{A_1}^{-1}{\mathbf Z}_{A_1}'\bolds\delta_n\|_2 &\le& 
\|{\mathbf C}_{A_1}^{-1} \|_2 \cdot\|n^{-1}{\mathbf Z}_{A_1}'{\mathbf
Z}_{A_1}\|_2^{1/2}
\cdot
\|n^{-1}\bolds\delta_n\|_2 \nonumber\\[-8pt]\\[-8pt]
& = & O_p(1) \rho_{n1}^{-1} \rho_{n2}^{1/2} q^{1/2} m_n^{-d}.\nonumber
\end{eqnarray}
By Lemma \ref{LemD}, the third term
%
%
\begin{equation}\quad
\label{Tc} \max_{j \in A_1} n^{-1}\lambda_{n2} \|{\mathbf T}_{nj}
{\mathbf C}_{A_1}^{-1}\bolds\nu_n\|_2 \le n \lambda_{n2} \rho
_{n1}^{-1} \|\bolds\nu
_n\|_2 =
O_p(1)\rho_{n1}^{-1} n^{-1} \lambda_{n2} h_n.
\end{equation}
Thus, (\ref{KKTaa}) follows from (\ref{Qnote}), (\ref{Ta})--(\ref{Tc}) and condition (B2a).
\end{pf}
\begin{lemma}
\label{LemF} Under conditions \textup{(B1)}, \textup{(B2)}, \textup{(A3)} and \textup{(A4)},
%
%
\begin{equation}
\label{KKTbb} \mathrm{P}\bigl(
\|{\mathbf Z}_{j}'({\mathbf Y}-{\mathbf Z}_{A_1}\widehat{\bolds\beta
}_{nA_1})\|_2 >
\lambda_{n2} w_{nj}/2, \exists j\notin A_1 \bigr) \rightarrow0.
\end{equation}
\end{lemma}
\begin{pf}
By (\ref{EEb}), we have
%
%
\begin{equation}
\label{ZZZa} {\mathbf Z}_j' ({\mathbf Y}-{\mathbf Z}_{A_1}\widehat
{\bolds\beta}_{nA_1} ) =
{\mathbf Z}_j'{\mathbf H}_n \bolds\varepsilon_n + {\mathbf
Z}_j'{\mathbf H}_n \bolds\delta_n+ \lambda n^{-1}
{\mathbf Z}_j'{\mathbf Z}_{A_1}{\mathbf C}_{A_1}^{-1}\bolds\nu_n.
\end{equation}
Recall $s_n=p-q$ is the number of zero components in the
model. By Lemma \ref{LemC},
%
%
\begin{equation}
\label{MaxTTa} \mathrm{E}\Bigl({\max_{j \notin A_1}}\| n^{-1/2}{\mathbf
Z}_j'{\mathbf H}_n
\bolds\varepsilon_n\|_2 \Bigr) \le O(1) \{ \log(s_n m_n)
\}^{1/2}.
\end{equation}
Since $w_{nj} = \|\widehat{\bolds\beta}_{nj}\|^{-1} = O_p(r_n)$ for
$j \notin
A_1$ and by (\ref{MaxTTa}), for the first term on the right-hand
side of (\ref{ZZZa}), we have
%
%
\begin{eqnarray}
\label{TTa}
&&
\mathrm{P}(\| {\mathbf Z}_j'{\mathbf H}_n \bolds
\varepsilon_n\|_2 >
\lambda_{n2} w_{nj}/6, \exists j \notin A_1 )
\nonumber\\
&&\qquad\le
\mathrm{P}(\| {\mathbf Z}_j'{\mathbf H}_n \bolds\varepsilon_n\|_2 >
C \lambda_{n2} r_n, \exists j
\notin
A_1 )+o(1) \nonumber\\[-8pt]\\[-8pt]
&&\qquad= \mathrm{P}\Bigl({\max_{j \notin A_1}}\| n^{-1/2}{\mathbf Z}_j'{\mathbf
H}_n \bolds\varepsilon_n\|_2
> C n^{-1/2}
\lambda_{n2} r_n \Bigr)+o(1) \nonumber\\
&&\qquad \le O(1) \frac{n^{1/2} \{ \log(s_n m_n)
\}^{1/2}}{C\lambda_{n2} r_n} +o(1).\nonumber
\end{eqnarray}
By (\ref{DeltaA}), the second term on the right-hand side of
(\ref{ZZZa})
%
%
\begin{eqnarray}
\label{TTb} {\max_{j \notin A_1}} \|{\mathbf Z}_j'{\mathbf H}_n \bolds
\delta_n\|_2&\le&
{n^{1/2} \max_{j \notin A_1}} \|n^{-1}{\mathbf Z}_j'{\mathbf Z}_j\|
_2^{1/2} \cdot
\|{\mathbf H}_n\|_2 \cdot\|\bolds\delta_n\|_2\nonumber\\[-8pt]\\[-8pt]
&=&O(1) n \rho_{n3}^{1/2}q^{1/2}
m_n^{-d}.\nonumber
\end{eqnarray}
By Lemma \ref{LemD}, the third term on the right-hand side of
(\ref{ZZZa})
%
%
\begin{eqnarray}
\label{TTc}
&&
\max_{j \notin A_1} \lambda_{n2} n^{-1} \|{\mathbf Z}_j
{\mathbf Z}_{A_1}{\mathbf C}_{A_1}^{-1} \bolds\nu_n\|_2 \nonumber\\
&&\qquad\le
{\lambda_{n2} \max_{j \in
A_1}}\|n^{-1/2}{\mathbf Z}_j\|_2 \cdot
\|n^{-1/2}{\mathbf Z}_{A_1}{\mathbf C}_{A_1}^{-1/2}\|_2\cdot\|{\mathbf
C}_{A_1}^{-1/2}\|_2
\cdot\|\bolds\nu_n\|_2 \\
&&\qquad= \lambda_{n2} \rho_{n3}^{1/2} \rho_{n1}^{-1/2} O_p(
qb_{n1}^{-1}).\nonumber
\end{eqnarray}
Therefore, (\ref{KKTbb}) follows from (\ref{Qnote}), (\ref{TTa}),
(\ref
{TTb}), (\ref{TTc}) and condition (B2b).
\end{pf}
\begin{pf*}{Proof of Theorem \ref{ThmD}}
The proof is similar to that
of Theorem \ref{ThmB} and is omitted.
\end{pf*}
\end{appendix}

\section*{Acknowledgments}
The authors wish to thank the Editor,
Associate Editor and two anonymous referees for their helpful
comments.

\printaddresses

\end{document}